\input amstex
\input pictex
\documentstyle{amsppt}
\magnification=\magstep1
\TagsOnRight
\topmatter

\hsize=6.5truein
\vsize=9truein

\def\bs{\baselineskip=.35truein}
\def\dist{{\text{Dist\,}}}
\def\dmax{{d_{\text max}\,}}
\def\rplus{{\Bbb R}^+}
\def\reals{{\Bbb R}}

\def\as {a{.}s{.}}
\def\iid{i{.}i{.}d{.}}
\def\rd {{\Bbb R}^d}
\def\zd {{\Bbb Z}^d}
\def\field{{\Cal F}}

\def\integers{{\Bbb Z}}

\def\prob{P}

\def\ex{E}

\def\var{{\text{Var\,}}}
\def\std{{\text{Std\,}}}

\def\qequals{\ = \ }
\def\qle{\ \le\ }

\def\qge{\ \ge\ }
\def\qg{\ >\ }
\def\e{\epsilon}

\def\x{{\ell \hat e_1}}
\def\d{\delta}
\def\D{\Delta}
\def\g{\gamma}
\def\G{\Gamma}
\def\a{\alpha}
\def\b{\beta}
\def\l{\lambda}

\def\s{\sigma}

\def\path{{$Q$-path}}

\def\qed{{$\ \blacksquare$}}

\def\eqab{ab}
\def\eqac{ac}
\def\eqad{ad}
\def\eqae{ae}
\def\eqaf{af}
\def\eqag{ag}
\def\eqah{ah}
\def\eqai{ai}
\def\eqaj{aj}
\def\eqak{ak}
\def\eqal{al}
\def\eqam{am}
\def\eqan{an}
\def\eqao{ao}
\def\eqap{ap}
\def\eqaq{aq}
\def\eqar{ar}
\def\eqas{as}
\def\eqat{at}
\def\eqau{au}
\def\eqav{av}
\def\eqaw{aw}
\def\eqax{ax}
\def\eqay{ay}

\def\eqbb{bb}
\def\eqbc{bc}
\def\eqbd{bd}
\def\eqbe{be}
\def\eqbf{bf}
\def\eqbg{bg}
\def\eqbh{bh}
\def\eqbi{bi}
\def\eqbj{bj}
\def\eqbk{bk}
\def\eqbl{bl}
\def\eqbm{bm}
\def\eqbn{bn}
\def\eqbo{bo}
\def\eqbp{bp}

\def\eqbr{br}

\def\eqbu{bu}

\def\eqbw{bw}
\def\eqbx{bx}

\def\eqab{3.1}
\def\eqac{3.3}
\def\eqad{3.4}
\def\eqae{3.5}
\def\eqaf{3.6}
\def\eqag{3.8}
\def\eqah{3.9}
\def\eqai{3.10}
\def\eqaj{3.11}
\def\eqak{3.19}
\def\eqal{3.20}
\def\eqam{3.21}
\def\eqan{3.22}
\def\eqao{3.12}
\def\eqap{3.13}
\def\eqaq{3.14}
\def\eqar{3.15}
\def\eqas{3.23}
\def\eqat{3.24}
\def\eqau{3.25}
\def\eqav{3.26}
\def\eqaw{3.16}
\def\eqax{3.17}
\def\eqay{3.18}

\def\eqbb{5.7}
\def\eqbc{5.8}
\def\eqbd{5.9}
\def\eqbe{5.10}
\def\eqbf{5.11}
\def\eqbg{5.12}
\def\eqbh{5.13}
\def\eqbi{5.14}
\def\eqbj{5.15}
\def\eqbk{5.16}
\def\eqbl{5.17}
\def\eqbm{5.18}
\def\eqbn{5.19}
\def\eqbo{5.20}
\def\eqbp{5.21}

\def\eqbr{3.29}

\def\eqbu{3.2}

\def\eqbw{3.27}
\def\eqbnew{3.28}
\def\eqbx{3.7}

\def\eqca{4.1}
\def\eqcb{4.2}
\def\eqcc{4.3}
\def\eqcd{4.4}
\def\eqce{4.5}
\def\eqcf{4.6}
\def\eqcg{4.7}
\def\eqch{4.8}
\def\eqci{4.9}

\def\cb{{C_1}}
\def\cs{{C_0}}

\def\kb{{\kappa_2}}  
\def\kc{{\kappa_3}}  
\def\kd{{\kappa_4}}  
\def\ke{{\kappa_5}}  
\def\kf{{\kappa_1}}  
\def\ka{{\kappa_1}}

\def\B{{\Cal B}}
\def\F{\tilde{F}}
\def\W{{\Cal W}}
\def\tR{\tilde{R}}
\def\tM{\tilde{M}}
\def\tD{\tilde{D}}
\def\M{M_{\alpha}}

\def\h{{f}}

\title    Geodesics and Spanning Trees for Euclidean First-Passage Percolation
\endtitle

\rightheadtext{Geodesics and Spanning Trees for Euclidean FPP}

\author   C. Douglas Howard$^1$ and Charles M. Newman$^2$  
\endauthor

\leftheadtext{C. Douglas Howard and Charles M. Newman}

\address  Baruch College,
          Box G0930
          17 Lexington Avenue
          New York, NY 10010, USA  
\endaddress

\email    dhoward\@baruch.cuny.edu   
\endemail

\address  Courant Institute of Mathematical Sciences,       
          New York University,
          251 Mercer Street,
          New York, NY 10012, USA   
\endaddress

\email    newman\@courant.nyu.edu    
\endemail  




\thanks $^1$  Research partially supported by NSF Grant DMS-98-15226.
\endthanks

\thanks $^2$  Research partially supported by NSF Grant DMS-98-03267.
\endthanks

\abstract The metric $D_\a(q,q')$ on the set $Q$ of 
particle locations of a homogeneous Poisson
process on $\rd$, defined as the infimum 
of $(\sum_i |q_i - q_{i+1}|^\a)^{1/\a}$ over sequences
in $Q$ starting with $q$ and ending with 
$q'$ (where $|\cdot|$ denotes Euclidean distance) has
nontrivial geodesics when $\a>1$. The cases $1<\a<\infty$ 
are the Euclidean first-passage percolation
(FPP) models introduced earlier by the authors 
while the geodesics in the case $\a=\infty$ are
exactly the paths from the Euclidean minimal spanning 
trees/forests of Aldous and Steele.
We compare and contrast results and conjectures 
for these two situations.
New results for $1<\a<\infty$ (and any $d$) include inequalities on the 
fluctuation exponents
for the metric ($\chi\le 1/2$) and for the geodesics 
($\xi \le 3/4$) in strong enough versions
to yield conclusions not yet obtained for lattice FPP:  
almost surely, every semi-infinite geodesic
has an asymptotic direction and every direction has 
a semi-infinite geodesic (from every $q$). For
$d=2$ and $2\le\a<\infty$, further results follow 
concerning spanning trees of semi-infinite geodesics
and related random surfaces.
\endabstract  

\endtopmatter

\document
\bs
\head{0.  Introduction}\endhead

There is an extensive literature (see [Y] for a survey) concerning combinatorial
optimization in which some functional based on the Euclidean distances $|q-q'|$ between
random points in $\rd$ is minimized.  Familiar examples include the total length in the
travelling salesman problem and in the minimal spanning tree.  In 
[HoN1], the authors introduced
another family of such functionals in order to 
obtain Euclidean versions of the
first-passage percolation (FPP) models originally defined in the context of the $\zd$ lattice by Hammersley 
and Welsh [HW]. (We remark that other Euclidean FPP models were 
introduced by Vahidi-Asl and Wierman [VW1, VW2] and studied by them and by
Serafini [Se].)
The focus of this paper is on these Euclidean FPP models from two perspectives.  First, we
survey a number of results and conjectures about 
these models with special emphasis on 
contrasts to the closely related but very different minimal spanning tree/forest of Aldous
and Steele [AS].  Second, we derive a number of new results about Euclidean FPP, and explain
why some of these go well beyond what has been proved for lattice FPP.  It is our hope that
the reader will find the pedagogical and research aspects of the paper to be complementary rather 
than antagonistic.

We define, for $r=(q_1,\dots,q_k)$ a finite 
sequence of points in $\rd$ and for $\a>0$
(usually we take $\a>1$),
$$
T_\a(r) = \sum_{j=1}^{k-1} |q_j-q_{j+1}|^\a,\ \  C_\a(r) = (T_\a(r))^{1/\a};
\eqno(0.1)
$$
for $\a=\infty$, we set
$$
C_\infty(r) = \max\{|q_j-q_{j+1}| : 1\le j < k\}.
\eqno(0.2)
$$
Starting from some (random) 
set of points $\tilde Q$ in $\rd$, we fix some $q$ and $q'$ in
$\tilde Q$ and then consider the combinatorial optimization problem of obtaining
$D_{\a,\tilde Q}(q,q') \equiv \inf\{C_\a(r)\}$ where the infimum is over all finite sequences
$r$ in $\tilde Q$ with $q_1 = q$ and $q_k = q'$ where $k$ is the (arbitrary) length of $r$.
When $\tilde Q$ is finite (e.g., $N$ independent uniformly distributed points in a cube of
volume $N$) there is of course some minimizing $r$ that yields the infimum, but we will be
interested in the case where $\tilde Q$ is a 
homogeneous Poisson point process on all of $\rd$ (corresponding
to $N\to\infty$) and then the issue of a minimizing $r$ is less trivial.

This issue is closely related to that of the existence of a geodesic path between $q$ and
$q'$ for the metric (when $\a\ge1$) $D_{\a,\tilde Q}$.  It turns out that the existence of
such a geodesic between arbitrary points $q$ and $q'$ is no problem for the Euclidean FPP
models where $1\le \a <\infty$, but for Euclidean minimal spanning trees, which as we shall
see correspond to $\a=\infty$, this 
is a serious issue which is not yet resolved for $d>2$.

In the next section of the paper, 
we give precise definitions of geodesics (finite and infinite), explain
why finite geodesics between arbitrary $q$ and $q'$
always exist when $\a<\infty$ and why they may not exist when $\a=\infty$.
We then review previous results for both lattice and Euclidean
FPP and state our new results
concerning the existence, nature and use of semi-infinite geodesics.
The latter are based on new estimates concerning 
the two exponents, $\chi$ and $\xi$, describing
respectively the fluctuation of the 
metric and of its geodesics.  These estimates are presented 
(and their relation to related results for lattice FPP is discussed) in
Section 2 and are used there to prove the new results of 
Section 1.  In
Sections 3 and 4,
the fluctuation exponent estimates are proved.
Some technical lemmas are given in Section 5.

\head{1. Geodesics and Spanning Trees}\endhead

Although our primary interest is in lattice and
Euclidean FPP and Euclidean minimal spanning trees/forests, we
will present the 
basic definitions in the general context 
of a countable set $\tilde Q$ (in our concrete examples,
this will be a subset of $\rd$ with $d\ge2$) 
and a function $\tau:\tilde Q\times \tilde Q
\to [0,\infty]$ (e.g., $\tau(q,q') = |q-q'|^\a$).  
We insist that $\tau(q,q)=0$ for every
$q\in\tilde Q$ and that $\tau(q,q')>0$ 
when $q\not=q'$ (although this latter condition
can be relaxed, e.g., in lattice FPP models).  
In our examples, $\tau(q,q') = \tau(q',q)$,
but this would not be so in directed (or oriented) FPP models.

A path $r$ is a sequence $(q_i:i\in I)$ that is indexed by an 
interval $I$ in $\integers$; it is
finite, semi-infinite or doubly-infinite according to the index set $I$. 
(For semi-infinite paths, we generally take $I$ infinite to the right, 
i.e., of the form $(i_0+1, i_0+2,\dots)$.) 
We also define a segment of a path
$r=(q_i:i\in I)$ to be any subpath
$r'=(q_i:i\in J)$ with $J$ a
sub-interval of $I$.  We call a path self-avoiding if
$q_i \not= q_j$ for any $i\not=j\in I$.

To each $i\in I$ (such that
$i+1\in I$) we associate $\tau_i = \tau(q_i,q_{i+1})$ and 
to each finite path
$r=(q_{i_0+1},\dots,q_{i_0+k})$ of length $k>1$, we associate a cost function
$\tilde C(r) = \tilde C_k(\tau(q_{i_0+1},q_{i_0+2}),\dots,
\tau(q_{i_0+k-1},q_{i_0+k}))\in[0,\infty]$
that is subadditive: for $k'\ge k>1$,
$$
\tilde C(\tau_1 ,\dots,\tau_{k-1},\tau_k ,\dots,\tau_{k'}) \qle 
   \tilde C(\tau_1,\dots,\tau_{k-1}) +  
   \tilde C(\tau_k ,\dots,\tau_{k'}).
\tag1.1
$$
(For a path $r$ of length 0, we take $\tilde C(r) = 0$.)  
Equivalently, in terms of a
path $r=(q,\dots,\hat q,\dots,q')$ from $q$ to $q'$ passing through $\hat q$ 
and thought of
as the concatenation of $r_1 = (q,\dots,\hat q)$ and $r_2=(\hat q,\dots,q')$, we have
$$
\tilde C(r) \qle \tilde C(r_1) + \tilde C(r_2).\tag1.2
$$
We also assume that $\tilde C (\tau_1,\dots,\tau_n) = \infty$ 
if and only if some $\tau_i=\infty$
and that\break
$\tilde C(\tau_1,\dots,\tau_n) = \tilde C(\tau_1,\dots,\tau_{j-1},\tau_{j+1},\dots,\tau_n)$
if $\tau_j=0$.  
The examples we consider are: 
$\tilde C = \sum_i\tau_i^{\,\a}$ or $(\sum_i\tau_i^{\,\a})^{1/\a}$
(with $\tau(q,q')=|q-q'|$);  $(\sum_i\tau_i)^{1/\a}$ 
(with $\tau(q,q')=|q-q'|^\a$); and $\max_i \tau_i$.
Taking $\a\ge1$ yields (1.1) and (1.2).

In the usual lattice FPP models (see, e.g., [Ke1]), $\tilde Q = \zd$ 
(or a subset of $\zd$) and 
$\tau(q,q') <\infty$ if and only if $q$ and $q'$ are nearest
neighbors on $\zd$; for such pairs, the $\tau(q,q')$'s are
i.i.d. random variables.  
In the Euclidean models of Vahidi-Asl and Wierman, $\tau(q,q') < \infty$
if and only if $q$ and $q'$ 
are neighboring points in the Voronoi or Delaunay graph associated
with $\tilde Q\subset \rd$.  
In our abstract setting one can define a graph $G$ with vertex set
$\tilde Q$ and edge set consisting 
of those $\{q,q'\}$ with $\tau(q,q') < \infty$. The assumption
(1.1) (or equivalently (1.2)) is important 
because it yields the triangle inequality for a natural metric
defined on each connected component of this graph as follows.

\proclaim{Definition}  Given $q,q'\in \tilde Q$, let $R(q,q')$ 
denote the set of all finite paths starting at $q$ and 
ending at $q'$ and define 
$$
\tilde D(q,q') \qequals \inf_{r\in R(q,q')} \tilde C(r) \tag1.3
$$
(or $\infty$, if $R(q,q')$ is empty).
\endproclaim

Note that $\tilde D(q,q)=0$ and that (1.2) yields 
the triangle inequality,
$$
\tilde D(q,q') \qle \tilde D(q,\hat q) + \tilde D(\hat q, q').
$$
In our abstract setting, $\tilde D$ may not 
be a metric (but only a pseudometric) 
if in taking the infimum 
in (1.3), $\tilde D(q,q')=0$ for some
$q\not=q'$.  This  happens, 
for example, with $\tilde C(r) = (\sum_i|q_i-q_{i+1}|^2)^{1/2}$
if $\tilde Q$ is dense in $\rd$.  $\tilde D$ will 
in fact be a metric in all of our
examples because 
$\tilde C_k(\tau_1,\dots,\tau_{k-1})\ge\tilde C_1(\tau_1)$ and
$$
\inf_{q'\not=q} \tilde C_1(\tau(q,q')) 
> 0 \text{\quad for all $q\in\tilde Q$}.\tag1.4
$$

\proclaim{Definition} A finite
path $r$ starting at $q$ and ending at $q'$ is said to be 
minimizing if the infimum in
(1.3) is finite and is achieved by $r$, i.e., if 
$\tilde D(q,q') = \tilde C(r) <\infty$.
A (finite, semi-infinite or doubly infinite) path
$r=(q_i:i\in I)$ is said to be a geodesic 
if it is self-avoiding and if every finite segment of $r$ is minimizing.
\endproclaim

We note that in all our examples 
{\it except\/} those with $\tilde C(r) = \max_i \tau_i$, every self-avoiding
minimizing finite path $r$ is automatically a geodesic.  
This is because if $r_{(1)}$ were a
non-minimizing segment of $r$, then representing $r$ as a concatenation of 
$r_{(0)}$, $r_{(1)}$ and $r_{(2)}$,
we could replace $r_{(1)}$ by an $r_{(1)}'$ 
(with the same endpoints as $r_{(1)}$ but with 
$\tilde C(r_{(1)}') < \tilde C(r_{(1))}$) 
and thus obtain an $r'$ with the same endpoints as $r$ and
with $\tilde C(r') < \tilde C(r)$ contradicting the minimizing property of $r$.

\subhead{1.1 Euclidean FPP}\endsubhead
In this subsection, 
we restrict attention to the Euclidean FPP models  of [HoN1] 
where $\tilde Q = Q$, the set of particle locations of a homogeneous
Poisson process of unit density on $\rd$, 
and $\tilde C(r) = C_\a(r) = (\sum_j |q_j-q_{j+1}|^\a)^{1/\a}$ with
$1\le \a < \infty$.  We 
denote the corresponding metric $\tilde D$ by $D_\a$.  Within our general
framework, one  may (i) set $\tau(q,q') = |q-q'|$ and 
$\tilde C(\tau_1,\dots,\tau_n) = (\sum_j \tau_j^{\,\a})^{1/\a}$,
or alternatively (ii)  set $\tau(q,q') = |q-q'|^\a$ and 
$\tilde C(\tau_1,\dots,\tau_n) = (\sum_j \tau_j)^{1/\a}$.
(Indeed, as far as the geodesics are concerned, one could instead (iii) set 
$\tau(q,q')= |q-q'|^\a$ and 
$\tilde C(r) = T_\a(r) = \sum_j |q_j - q_{j+1}|^\a$
so that $\tilde C (\tau_1,\dots,\tau_n) = \sum_j \tau_j$.  
The latter is best for 
comparing Euclidean FPP with 
lattice FPP while (i) is best for letting $\a\to\infty$ so that
both $\tau(q,q')$ and $\tilde C(r)$ have a limit.)

When $\a = 1$ (or in version (iii) above
when $0<\a\le1$) and $d\ge2$, since (almost surely)
no three points of $Q$ are collinear, 
it follows that for any distinct $q,q'\in Q$, the
unique geodesic between them is 
the trivial one going from $q$ to $q'$ in one step --- i.e.,
$r=(q,q')$.  To get nontrivial 
geodesics we need $\a>1$.  The next proposition states that
for $\a\not=\infty$, geodesics exist 
between all pairs of points (and are unique).  It is 
the Euclidean analog of a standard result in lattice FPP 
(see, e.g., [SmW])
with a similar
proof, which we sketch.  Our focus for $\alpha < \infty$ 
will then be on the asymptotic behavior of the finite 
geodesic between $q$ and $q'$ as $|q-q'|\to\infty$ 
and on the existence, nature and abundance of infinite
geodesics.  As we shall see 
in the next subsection, when $\a=\infty$, even the existence of
finite geodesics is nontrivial.

\proclaim{Proposition 1.1}  In Euclidean FPP 
with $d\ge2$ and $1\le\a <\infty$, there is almost surely
a unique geodesic $\M(q,q')$ between every pair of distinct points $q,q'\in Q$.
\endproclaim

\demo{Proof} The uniqueness follows because if $r$ 
and $r'$ were different self-avoiding paths from
$q$ to $q'$ with $C_\a(r) = C_\a(r')$, there would be two disjoint sets
$\{\{q_i,\bar q_i\}:i=1,\dots,m;\,q_i \neq \bar q_i \}$ and  
$\{\{q'_j,\bar q'_j\}:j=1,\dots,m';\,q'_j \neq \bar q'_j\}$
with
$\sum_i |\bar q_i - q_i|^\a = \sum_j |\bar q'_j - q'_j|^\a$.
But that occurs with zero probability.

To prove existence, note that the 
intersection of $Q$ with the Euclidean ball
${\B}(0,K) \equiv \{x\in \rd : |x|\le K\}$ 
is, for any $K<\infty$, almost surely finite.  
We define for $q\in Q$,
$$
d_\a(q,K) \qequals \inf\{D_\a(q,q') : q'\in Q \setminus {\B}(0,K)\}\tag1.5
$$
and claim that for every $q\in Q$,
$$
\lim_{K\to\infty} D_\a(q,K) \qequals \infty\ \ \text{\as}\tag1.6
$$	
This implies that 
for given $q$, $q'$ and then some sufficiently large $K$, any
$r$ from $q$ to $q'$ that exits $\B(0,K)$ 
has $C_\a(r) > |q-q'|$ 
and hence the infimum in (1.3) must be achieved within the 
{\it finite\/} collection of (self-avoiding) paths 
staying in
$\B(0,K)$.  The claim (1.6) is proved by appeal to a standard 
(continuum) percolation result (see [ZS1, ZS2]) --- namely 
that for some sufficiently small $\e>0$, any semi-infinite self-avoiding
path in $Q$ must make infinitely 
many steps with $|q-q'|>\e$.  This easily yields (1.6). \qed
\enddemo

As we shall see, analyzing the existence 
and nature of infinite geodesics can be 
difficult.  
But the following proposition, which shows that there
is at least one semi-infinite geodesic 
starting from each $q\in Q$, is not hard.

\proclaim{Proposition 1.2}  Suppose $d\ge2$ 
and $1<\a<\infty$.  For each $q\in Q$ define $R_\a(q)$ to
be the graph with vertex set $Q$ and edge set 
$\cup_{q'\in Q} \M(q,q')$.
Almost surely, for every $q\in Q$, $R_\a(q)$ 
is a spanning tree on $Q$ with every vertex
having finite degree; thus there is at least 
one semi-infinite geodesic starting from every $q$.
\endproclaim

\demo{Proof}  To see that $R_\a(q)$ is a spanning tree, order
$Q = \{q_1, q_2,\dots\}$ and note that inductively, for each $n$,
$\cup_{i=1}^n \M(q, q_i)$ is a tree because of the uniqueness
part of Proposition 1.1.  To justify the finite degree claim
(which we note is not valid when $\a=1$), it suffices to show
that for each $\tilde q\in Q$, there are \as\ only finitely
many $\bar q\in Q$ such that the single 
step path $r=(\tilde q,\bar q)$ is a geodesic.
This is a consequence of
$$
\lim_{K\to\infty} P[\text{$(\tilde q,\bar q)$ is a geodesic for some
$\tilde q,\bar q$ with $|\tilde q|\le 1$, $|\bar q|\ge K$}]\qequals 0,
\tag1.7 
$$
which itself follows from Lemma 5.2 (see (5.5)). We note that
the key geometric idea here is to define
$$
\W(a,b) \qequals \{x\in\rd : |a-x|^\a + |x-b|^\a \le |a-b|^\a\}\tag1.8
$$
and realize that
$(\tilde q, \bar q)$ cannot be a geodesic unless
$\W^\circ(\tilde q, \bar q)$, the interior of $\W(\tilde q, \bar q)$, is devoid
of Poisson particles. \qed 
\enddemo

When $\a=\infty$, there will also be at
least one semi-infinite geodesic starting from every $q$.
In that case, however, it is believed to be unique
(see Conjecture 1 below)--- unlike
when $\a<\infty$ as we discuss later.
A question apparently first posed (for lattice FPP) by 
H. Furstenberg (see p. 258 of [Ke1]) is: 
What about doubly infinite geodesics?  Here it is believed that
\as\ these do not exist both for $\a<\infty$ and $\a=\infty$.  We
shall see later the extent to which this has been proved.

\subhead{1.2 Minimal Spanning Trees and Forests}\endsubhead
In this subsection (and the rest of the paper)
we continue to take $\tilde Q=Q$, a homogeneous
Poisson process of unit density (except as noted) on $\rd$, but 
for now we take 
$\tilde C(r)=C_\infty(r) = \max_j |q_{j+1}-q_j|$ with the corresponding
metric $\tilde D(q,q') = D_\infty(q,q')$, the minimax of $|q_{j+1}-q_j|$
along paths $r\in R(q,q')$.

Let us denote by $R^*(q,\bar q)$ the set of all paths in $R(q,\bar q)$
that do not use the edge $\{q,\bar{q}\}$.  In order that the edge
$\{q,\bar q\}$ belong to some geodesic, it is necessary and sufficient
that $r=(q,\bar q)$ is itself a geodesic and, \as, this is true
if and only if 
$$
\text{$|q-\bar q| < C_\infty(r)$ for every $r\in R^*(q,\bar q)$.}
\tag1.9
$$
Following Alexander [Al2],
we make the following

\proclaim{Definition}  
$R_\infty$ is the graph with vertex set $Q$
and edge set consisting of those $\{q,\bar q\}$'s satisfying (1.9).
\endproclaim

The graph $R_\infty$ can have no loops because 
on any loop, the edge $\{q,\bar q\}$
with maximum $|q-\bar q|$ does not satisfy (1.9).
Thus $R_\infty$ is a forest (a union of one or more disjoint trees)
and contains at most one path between any $q,q'$.
Every finite geodesic $M_{\infty}(q,q')$ must be a path in $R_\infty$
and it is also not hard to see that every path in $R_\infty$
is a geodesic.  Thus, as in the $\a<\infty$ case, if a geodesic
$M_{\infty}(q,q')$ exists between $q$ and $q'$ it will be unique;
however, when $\a=\infty$, it may not exist.
$M_{\infty}(q,q')$ exists if and only if $q$ and $q'$
are in the same connected component of $R_\infty$ (if they
are in different components, we set $M_{\infty}(q,q')=\emptyset$).  
Thus geodesics
exist between every pair $q$, $q'$ in $Q$ if and only if
$R_\infty$ is a single (spanning) tree.  

\proclaim{Definition}  $R_\infty(q)$ is the graph with vertex set $Q$
and edge set $\cup_{q' \in Q} M_{\infty}(q,q')$.
\endproclaim

Clearly $R_\infty(q)$ is just 
the connected component of $q$ in $R_\infty$;  it is not
hard to see that each $R_\infty(q)$ must be an infinite tree.
If $R_\infty$ is a single tree, then (unlike when $\a<\infty$)
the $R_\infty(q)$'s are all the {\it same\/} spanning tree.

It is shown in [Al2] that $R_\infty$ is the same 
as the minimal spanning forest (MSF) constructed
by Aldous and Steele [AS] as follows: For
$K<\infty$, let $R_\infty^K$ denote the spanning tree of
$Q\cap \B(0,K)$ that minimizes the sum of $|q-\bar q|$ over
all edges $\{q,\bar q\}$ in the tree; then $R^K_\infty\to R_\infty$
as $K\to\infty$.  There are two obvious qualitative issues
concerning $R_\infty$. Is it a single spanning tree or not?
How many 
different semi-infinite
geodesics start from $q$?  This number, which is clearly the
same for all $q$'s in any fixed connected component of
$R_\infty$, equals the number of (topological) ends of
the component. (An end is an equivalence class of semi-infinite
paths in $R_\infty$ that agree except for finite initial segments.)
As to the number of ends, Alexander's
results [Al2] combined with a natural conjecture 
about continuum percolation lead to the following
(the natural conjecture is that at the critical radius $R_c^*$ for 
overlapping balls of fixed radius centered at points of $Q$
to form infinite clusters, there \as\ is no infinite cluster; for an
extensive presentation of rigorous 
results about continuum percolation, see [MR]):

$\bullet$\, {\bf{Conjecture 1 [Al2].}} For any $d\ge2$, $R_\infty$
contains exactly one semi-infinite geodesic from each $q$.

Note that this includes the conjecture that there are
no doubly infinite geodesics.  The latter conjecture will persist
for $\a<\infty$ even though Conjecture 1 will not.
As to the other issue, the natural extension from the lattice case
of a conjecture of Newman and Stein [NewS1, NewS2] is

$\bullet$\, {\bf{Conjecture 2 [NewS1, NewS2].}} For $d<8$ 
(and perhaps also $d=8$), $R_\infty$
is a single spanning tree;  for $d>8$, $R_\infty$ has (infinitely)
many connected components.

The only dimension where these conjectures have been verified
is $d=2$, as stated in the next theorem. However we note that
Conjecture 1 has been verified in lattice
models also for large $d$---see Example 2.7 of [Al2].
For general $d$, it has been proved [Al2]
that at most one component of $R_\infty$ has two ends
and all others have a single end.

\proclaim{Theorem 1.3 [AlM, Al2]}  For
$d=2$, $R_\infty$ is a single spanning tree with one 
end.
\endproclaim

In the next two subsections, we investigate the quite different qualitative
nature of semi-infinite geodesics when $\a<\infty$.  There will be
many more infinite geodesics from each $q$ and they will be asymptotically
fairly regular.  The irregularity of the infinite (or very long)
paths in $R_\infty$ is itself an interesting object of study.
One way to pursue this issue is to consider for each $x$
in $\rd$ the (unique for $d=2$ or
under Conjecture 1) infinite
path in $R_\infty$ starting from (the $q$ closest to)
$x$, in the model with Poisson density $1/\delta^d$, as a random
curve in $\rd$ and study its subsequence 
limits in distribution as $\delta\to0$.
Some interesting results in this regard (especially for $d=2$) have
been obtained in [ABNW] using technical methods
from [AB]. There are also interesting results on such scaling limits
for other random spanning tree models in [ABNW] and [S].

\subhead{1.3  Previous Results for Euclidean FPP}\endsubhead
There are two types of previously known results.  The first, valid
for all $d$ and $1<\a<\infty$, concerns the asymptotic shape of large
balls based on the metric $D_\a$.  The second, proved only for $d=2$ 
and $2\le\a <\infty$, concerns semi-infinite geodesics 
$r=(q_1 , q_2 ,\dots)$ with a specified asymptotic direction
$\hat x$ --- i.e., such that $q_k/|q_k|\to\hat x \in S^{d-1}$ as $k\to\infty$.
We will call such an $r$ an $\hat x$-geodesic.  Doubly infinite geodesics
$(\dots,q_{-1},q_0,q_1,\dots)$ such that $q_k/|q_k|\to\hat x$ 
(resp., $\hat y$) as $k\to\infty$ (resp., $-\infty$) will be
called $(\hat x,\hat y)$-geodesics.

Both types of results were originally derived in [HoN1] as analogs
of corresponding lattice FPP results.  The first type differs
from the lattice case in that the asymptotic shape is exactly
a Euclidean ball (because of the statistical Euclidean invariance
of the homogeneous Poisson point process).  The significance of this 
difference for our new results will be discussed in Section 2 of this
paper. We present the shape theorem result in a slightly different
form than the one of [HoN1]; in Section 2 (Theorem 2.3)
we improve this result.  For $x\in\rd$, denote by $q(x)$ the
Poisson particle location in $Q$ closest to $x$ (with any fixed
rule for breaking ties).  Then for $s>0$, let
$B_\a(x,s)\equiv\{q'\in Q:D_\a(q(x),q') \le s\}$ denote the
ball in $Q$ of radius $s$ centered at $q(x)$, using the metric
$D_\a$.

\proclaim{Theorem 1.4 [HoN1]} For any $\a\in(1,\infty)$ and $d\ge 2$,
there exists  $\mu\in(0,\infty)$ depending on $\a$ and $d$, such that 
with ${\Cal B}_0 \equiv {\Cal B}(0,\mu^{-1})$
the following is true almost surely.  For any $\e\in(0,1)$,
$$
Q\cap (1-\e)s{\Cal B}_0 \ \subset\ B_\a(0,s^{1/\a})\ \subset
\ (1+\e)s{\Cal B}_0
\tag1.10
$$
for all sufficiently large $s$.
\endproclaim


There are many natural questions one
can ask about semi-infinite geodesics.  We may focus on some $q\in Q$
(e.g., $q(0)$, the particle nearest the origin) and consider
(for a fixed $\a$), the set $G_\a(q)$ of semi-infinite geodesics starting
from $q$.  $G_\a(q)$ is of course just the set of 
semi-infinite paths starting from $q$ in the spanning tree $R_\a(q)$ defined in
Subsection 1.1 above, so that (for $1<\a<\infty$, according
to Proposition 1.2) $G_\a(q)$ is nonempty.

When  $\a=\infty$, as discussed in Subsection 1.2, $R_\a(q)$ may
not be spanning (for large enough $d$), but it is still
an infinite tree of finite degree at each
vertex, so $G_\infty(q)$ is also nonempty.
For $\a=\infty$ and $d=2$, according to Theorem 1.3, $G_\infty(q)$
consists of a single infinite geodesic and further for any
$q$ and $q'$, the (unique) semi-infinite geodesics $r$ and $r'$
starting from $q$ and $q'$ coalesce; i.e., there is a unique 
$\bar q\in Q$ (which may be $q$ or $q'$) such that $r$ (resp., $r'$)
is the concatenation of a path $\tilde r$ from $q$ to $\bar q$
(resp., $\tilde r'$ from $q'$ to $\bar q$) with the semi-infinite
geodesic $\bar r$ starting from $\bar q$, while $\tilde r$ and $\tilde r'$
are disjoint except for $\bar q$.  It is not hard to show (using the
statistical rotational invariance of the Poisson point process) that
here the semi-infinite geodesics cannot have an asymptotic direction
(indeed, that the set of subsequence limit points of $q_k/|q_k|$
along a semi-infinite geodesic must \as\ be all of the
unit circle).

But for $1<\a<\infty$ and arbitrary $d$, one expects rather different
answers to the following questions.

$\bullet$\,{\bf Question 1.}  Does every 
semi-infinite geodesic have an asymptotic
direction?

$\bullet$\,{\bf Question 2.}  Is $R_\a(q,\hat x)$, 
the set of $\hat x$-geodesics
starting from $q$, nonempty for every unit vector $\hat x$?

$\bullet$\,{\bf Question 3.} Are there some $\hat x$'s with more than one
$\hat x$-geodesic from some $q$ (i.e., with $R_\a(q,\hat x)$
bigger that a singleton)?

If the answer to Question 3 turns out to be ``Yes'', one may ask
a related but different question, whose answer could still be ``No'',
as follows.

$\bullet$\,{\bf Question 4.}  For a 
deterministic $\hat x$, can there be more than
one $\hat x$-geodesic from some $q$ and can there be
non-coalescing $\hat x$ geodesics from different $q$'s?

$\bullet$\,{\bf Question 5.}  Do doubly-infinite geodesics exist?

``Yes'' answers to Questions 
1, 2, and 3 are among the main new results of this paper
and will be stated as 
theorems in the following subsection.  Analogous results
for lattice FPP are still open problems (see [New1] and Section 2
of this paper).
The answer ``No'' to the fourth question was 
previously known, but only for restricted $d$ and $\a$
(it remains an open problem in general), as follows; 
the restriction on $\a$ will be discussed below:

\proclaim {Theorem 1.5 [HoN1]} For $d=2$, $2\le\a<\infty$ 
and every deterministic $\hat x$, 
\as\ there is no more than one $\hat x$ geodesic from any $q$ and \as\
any pair of $\hat x$ geodesics from  distinct $q,q'$ must coalesce. 
\endproclaim

We remark that there are lattice FPP 
analogs to this theorem (and the next), but these
have not been proved for {\it every\/} $\hat x$ [LN]; 
the best such result is due
to Zerner (Theorem 1.5 in [New2]).  As a consequence of the last theorem,
there was a partial answer 
to Question 5, stated as the next theorem.
The natural conjecture is that the correct answer to Question 5
is ``No'', certainly for $d=2$ and perhaps for all $d$. 
(See Chap. 1 of [New2] for a discussion of this conjecture
for lattice FPP and its equivalence (when $d=2$)
to nonexistence of nonconstant ground states for 
disordered Ising ferromagnets.
Other results in the lattice context are in [W].)

\proclaim {Theorem 1.6 [HoN1]} For $d=2$, $2\le\a<\infty$, and every deterministic $\hat x$ and $\hat y$,
\as\ there are no $(\hat x, \hat y)$-geodesics.
\endproclaim

An improvement of Theorem 1.6 (see Theorem 1.11) will be 
given below,
basically as a consequence of our answer to 
Question 1, but this improvement falls well
short of the conjecture that doubly infinite geodesics \as\ do not exist.

Behind the restriction to $\a\ge2$ 
in Theorem 1.5 and 1.6 is 
the following lemma (Lemma 5 of [HoN1]), which we will use later.

\proclaim{Lemma 1.7 [HoN1]} Suppose $r = (\dots,q_i,q_{i+1},\dots)$ 
and $r' = (\dots,q'_j,q'_{j+1},\dots)$
are two finite or infinite 
geodesics such that the closed line segments $\overline{q_i q_{i+1}}$ and
$\overline{q'_j q'_{j+1}}$ intersect.  If $d=2$ and $2\le\a<\infty$, then
$\{q_i,q_{i+1}\}$ and $\{q'_j,q'_{j+1}\}$ have at least one point in common.
\endproclaim

\subhead{1.4  New Results on Infinite 
Geodesics for Euclidean FPP}\endsubhead  The next three
theorems, among the main 
new results of this paper, are consequences of fluctuation
theorems presented in Section 2.  
The fluctuation theorems are of interest in their own right. 

\proclaim{Theorem 1.8} For $d\ge2$, 
and $1<\a<\infty$, \as: every semi-infinite geodesic
has an asymptotic direction.
\endproclaim

\proclaim{Theorem 1.9} For $d\ge2$ 
and $1<\a<\infty$, \as: for every $q\in Q$ and every
unit vector $\hat x$, there 
is at least one $\hat x$-geodesic starting from $q$.
\endproclaim

\proclaim{Theorem 1.10} For $d\ge2$ 
and $1<\a<\infty$, \as: for every $q\in Q$, the set $V(q)$ of
unit vectors $\hat x$ such 
that there is more than one $\hat x$-geodesic starting at
$q$ is dense in the unit sphere.
\endproclaim

{\it Remark.}  for $d=2$, is is not hard to show (by arguments like those 
used to prove Theorem 0 of [LN]) that
\as\ $V(q)$ is countable.  In general, 
whenever the answer to the first part of Question 4 is ``No'', then
by an application of Fubini's Theorem, 
the Lebesgue measure (on the unit sphere $S^{d-1}$) of $V(q)$ is
zero.  But the proof of 
Theorem 1.10 also shows 
that $V(q)$ must have Hausdorff dimension at least $d-2$.

Theorem 1.8 implies that every doubly-infinite geodesic must be an 
$(\hat x, \hat y)$-geodesic for some
$\hat x, \hat y\in S^{d-1}$.  But the proof of Theorem 1.8 implies a bit more.  We state this in the 
next theorem in combination with the result of Theorem 1.6.

\proclaim{Theorem 1.11}  For $d\ge2$ and $1<\a<\infty$, \as\ 
doubly infinite geodesics other than
$(\hat x, - \hat x)$-geodesics 
do not exist.  In addition, for $d=2$ and $2\le\a<\infty$, and any
deterministically chosen 
$\hat x$, \as\ $(\hat x,-\hat x)$-geodesics do not exist.
\endproclaim

Theorem 1.11 is a step in the direction of verifying 
the conjecture that, \as, doubly infinite geodesics do not exist.  
However, even for $d=2$ and
$2\le\a<\infty$, it does not 
prove the conjecture since it leaves open the possible existence
of $(\hat x, -\hat x)$-geodesics with $\hat x$ dependent on the
realization of $Q$.

In the rest of this subsection, we restrict 
attention to $2\le\a<\infty$ and $d=2$ and
explore some consequences of combining Theorems 1.5-1.11.  This is in the spirit of [New1] 
(see Theorem 1.1 of that reference and the preceding discussion there), where
the same issues were addressed, but only partially resolved,
in the lattice FPP context.

When $d=2$, $2\le\a<\infty$, $q\in Q$, and $\hat x$ is a deterministic
unit vector (in $S^1$), by Theorems 1.5 and 1.9, there \as\ exists
a unique $\hat x$-geodesic starting from $q$.  We denote this 
semi-infinite geodesic by $s_q(\hat x)$.  In analogy with
$R_\a(q)$, as defined in Proposition 1.2, (but with $q$ replaced
by ``a point at infinity reached in the direction $\hat x$'')
we define $R_\a(\hat x)$
to be the graph with vertex set $Q$ and every edge contained in
$\cup_{q'\in Q} s_{q'}(\hat x)$.  It follows from Theorem 1.5 that (\as)
$R_\a(\hat x)$ is a spanning tree on $Q$ (the coalescing part of Theorem
1.5 ensures that $R_\a(\hat x)$ has a single connected component).
Since every edge in $R_\a(\hat x)$ touching $q$ is part of some geodesic,
these edges belong to $R_\a(q)$ and hence, by Proposition 1.2, each
vertex in $R_\a(\hat x)$ has finite degree.  We combine these facts with
a few others in the following.

\proclaim{Theorem 1.12}  Suppose $d=2$, $2\le\a<\infty$, and
$\hat x$ is a deterministic unit vector (in $S^1$).  Then the
following are all valid \as. For
any $q\in Q$ and any $\bar q_1, \bar q_2,\dots\in Q$ such that
$\bar q_k / |\bar q_k| \to \hat x$, the finite geodesic $\M(q, \bar q_k)$
converges as $k\to\infty$ to the unique
$\hat x$-geodesic $s_q(\hat x)$ starting from $q$.  Thus the spanning
trees $R_\a(\bar q_k) \to R_\a(\hat x)$ as $k\to\infty$ (where the
edges of $R_\a(\hat x)$, as defined above,
are those in $\cup_{q\in Q} s_q(\hat x)$).
$R_\a(\hat x)$ is a spanning tree on $Q$ (with every vertex having
finite degree and) with a single infinite path from each $q$
(namely $s_q(\hat x)$);  $R_\a(\hat x)$ thus has a single
topological end.
\endproclaim

\demo{Proof of Theorem 1.12} The things that remain to be proved
are that $\M(q,\bar q_k)\to s_q(\hat x)$ and that $R_\a(\hat x)$
contains no infinite path from $q$ other than $s_q(\hat x)$. For
a small $\e>0$, let $\hat x_+(\e)$ (resp., $\hat x_-(\e)$)
be the unit vector obtained by rotating $\hat x$ by an angle
$\e$ in the clockwise (resp., counterclockwise) direction.
By Theorems 1.5 and 1.9, there \as\ exist unique semi-infinite 
geodesics $s_q(\hat x_{\pm}(\e))$ starting from $q$.  For a 
path $r=(q_1 , q_2 ,\dots)$ let us denote by $\bar r$
the union of the line segments $\overline{q_i\,q_{i+1}}$
(as a subset of ${\Bbb R}^2$). The paths $s_q(\hat x_+(\e))$ and
$s_q(\hat x_-(\e))$ bifurcate at some $q'$ (perhaps equal to $q$)
and then, by uniqueness of finite geodesics, have no further $Q$-particles
in common.  By Lemma 1.6, the sets $\overline{s_q(\hat x_+(\e))}$
and $\overline{s_q(\hat x_-(\e))}$ bifurcate at $q'$ and have no
further ${\Bbb R}^2$-points in common. Thus
${\Bbb R}^2 \setminus 
\{(\overline{s_q(\hat x_+(\e))}\cup\overline{s_q(\hat x_-(\e))}\}$ consists
of two connected components (one ``inside'' and one ``outside'') that
we will denote by $S_q^{\text{in}}(\hat x,\e)$ and
$S_q^{\text{out}}(\hat x,\e)$. The inside (resp., outside)
component is characterized by containing sequences $x_1, x_2,\dots$
in ${\Bbb R}^2$ such that $|x_j|\to\infty$ while the 
angle between $x_j/|x_j|$ and $\hat x$ converges to a point in 
$(-\e,\e)$ (resp., to a point outside $[-\e,\e]$).

Now, by Lemma 1.6 again (and the uniqueness of finite geodesics)
once $k$ is large enough that $\bar q_k\in S_q^{\text{in}}(\hat x,\e)$,
$\overline{\M(q,\bar q_k)}$ (except for its initial portion 
from $q$ to $q'$) must be entirely within the
closure of $S_q^{\text{in}}(\hat x,\e)$ and thus the same
must be true for any (subsequence) limit $\tilde r$ of $\M(q,\bar q_k)$.
Since this is true for every $\e>0$, it follows that such an $\tilde r$
(which is automatically a geodesic starting from $q$) must be an
$\hat x$-geodesic. But then by Theorem 1.5, $\tilde r$ is
\as\ $s_q(\hat x)$ as claimed.

Next suppose that $\hat r$ is an infinite path in $R_\a(\hat x)$ 
starting from $q$  and different than $s_q(\hat x)$.  We show that this leads
to a contradiction.  The path $\hat r$ must bifurcate from
$s_q(\hat x)$ at some $q'$ (possibly with $q'=q$) with 
no further $Q$-particles in common.  For any $q''$ on $\hat r$ after $q'$,
the concatenation of the segment of $\hat r$
from $q''$ to $q'$ and the infinite segment of $s_q(\hat x)$
starting at $q'$ (which is just $s_{q'}(\hat x)$) must be 
$s_{q''}(\hat x)$ since $s_{q''}(\hat x)$ and $s_q(\hat x)$
must coalesce somewhere and if it were not at $q'$, $R_\a(\hat x)$
would contain a loop. 
Let $q''_k$ denote an infinite sequence
of distinct such $q''$'s from $\hat r$
and let $r''$ be a limit of 
$s_{q''_k}(\hat x)$ --- which must exist since each $s_{q''_k}(\hat x)$
passes through $q'$ and contains $s_{q'}(\hat x)$.
Then $r''$ is a doubly infinite geodesic containing $s_{q'}(\hat x)$
and thus by the first part of Theorem 1.11 is an
$(\hat x, -\hat x)$-geodesic.  But this contradicts the second
half of Theorem 1.11, which completes the proof.\qed
\enddemo

Now that we have constructed in Theorem 1.12
the spanning tree $R_\a(\hat x)$ composed of the
$\hat x$-geodesics $s_q(\hat x)$, we may ask: what is it good for?
Following [New1], it can be used to study the surface 
of large balls in the metric space $(Q, D_\a)$ by means of certain (random)
``height functions'' on $Q$ (or on $\rd$).  For a fixed
$\a<\infty$ we replace $D_\a$ by the pseudometric on $\rd$,
$T_\a(x,y)\equiv D_\a(q(x),q(y))^\a$ (where $q(x)$ is the
closest $q\in Q$ to $x$) and look at the pseudometric balls, 
$\tilde B_\a(x,s) \equiv \{y\in \rd: T_\a(x,y)\le s\}$.
These are unions of Voronoi regions and are related to the balls
$B_\a(x,s)$ for the metric $D_\a$ (defined just above Theorem 1.4) by
$B_\a(x,s^{1/\a}) = \tilde B_\a(x,s)\cap Q$.  

What does a
large-radius ball $\tilde B_\a(x,s)$ look like when ``viewed from
its surface?''  A natural interpretation of this question, that
places the surface near the origin, is to consider the limit\
of $\tilde B_\a(\bar q_k, T_\a(\bar q_k,0))$ as $|\bar q_k|\to\infty$
with $\bar q_k/|\bar q_k|\to\hat x$.  Theorem 1.12 allows us
to analyze this limit in terms of a function $H^{\hat x}(q,q')$
on $Q\times Q$ defined as follows.  For $q, q'\in Q$, define
$W_{\hat x}(q,q')$ as the unique $q''$ in $Q$ where $s_{\hat x}(q)$
and $s_{\hat x}(q')$ coalesce ($W_{\hat x}$ might be $q$ or $q'$)
so that the path in $R_\a(\hat x)$ between $q$ and $q'$ is
the concatenation of $\M(q,W_{\hat x}(q,q'))$ and $\M(W_{\hat x}(q,q'),q')$.
The following is mostly a consequence of Theorem 1.12.

\proclaim{Theorem 1.13}  Suppose $d=2$, $2\le\a<\infty$, and
$\hat x$ is a deterministic direction.  Then the following
are valid \as:  For all $q,q'\in Q$ and any $\bar q_1, \bar q_2,\dots\in Q$
such that $\bar q_k/|\bar q_k|\to \hat x$, 
$$
H^{\hat x}(q,q') \ \equiv\ 
\lim_{k\to\infty} [T_\a(q,\bar q_k) - T_\a(q', \bar q_k)]
$$
exists and equals $T_\a(q,W_{\hat x}(q,q')) - T_\a(q',W_{\hat x}(q,q'))$.
The balls $\tilde B_\a(\bar q_k, T_\a(\bar q_k,0))$ converge as $k\to\infty$
to $\{y\in\rd:H^{\hat x}(q(y),q(0))\le 0\}$.  Furthermore,
$H^{\hat x}(\cdot,q(0))$ as a function on $Q$ satisfies
$$
H^{\hat x}(q,q(0)) \qequals \inf_{q'\not=q}[|q-q'|^\a + H^{\hat x}(q',q(0))]
\tag1.11
$$
and more generally for $Q_0$ any finite subset of $Q$ containing $q$
$$
H^{\hat x}(q,q(0)) \qequals \inf_{q'\in Q\setminus Q_0}[T_\a(q,q') + 
H^{\hat x}(q',q(0))]
\tag1.12
$$
\endproclaim

\demo{Proof}  The only claims that require any explanation are 
(1.11) and (1.12).  To prove (1.11), we let $q''$ denote the first
particle after $q$ on $s_{\hat x}(q)$ and note that by Theorem 1.12,
$$
\align
H^{\hat x}(q,q(0)) 
&\qequals \lim_{k\to\infty}[|q-q''|^\a + 
T_\a(q'',\bar q_k) - T_\a(q(0),\bar q_k)]\\
&\qequals |q-q''|^\a + H^{\hat x}(q'',q(0)).\tag1.13
\endalign
$$
This bounds $H^{\hat x}(q,q(0))$ below by the right side
of (1.11).  The opposite inequality easily follows from
$T_\a(q,\bar q_k) \ge \inf_{q'\not= q} (|q-q'|^\a + T_\a(q',\bar q_k))$.
The identity (1.12) is derived by quite similar arguments to those 
used for (1.11).\qed
\enddemo

We now consider the random field $H^{\hat x}(q(y),q(0))$.  It is
clear, at least on a heuristic level, that the asymptotic behavior
of its mean, as $|y|\to\infty$, is 
$-\mu(\a,2) (\hat x\cdot y)$ to leading order,
where $\mu(\a,d)$ is the inverse of the 
radius appearing in Theorem 1.4 and $\hat x\cdot y$ denotes the
standard Euclidean inner product.  When $\hat x\cdot y \not= 0$, it seems
reasonable that the variance of $H^{\hat x}(q(y),q(0))$ should have a 
leading order behavior similar to that of $T_\a(0,y)$ --- namely like
$|y|^{2\chi}$ (with $\chi=1/3$ conjectured for $d=2$, 
as discussed in Section 2). For
$\hat x\cdot y = 0$, where by symmetry $E[H^{\hat x}(q(y),q(0)] = 0$,
it seems that for $d=2$, one should expect the variance to grow 
faster, namely linearly in
$|y|$, and correspondingly the boundary of the region where 
$H^{\hat x}(q(x),q(0)) \le 0$ should fluctuate from the straight line
$y=t\hat y_0$ (where $\hat y_0 \cdot \hat x = 0$) by a distance
of order $\sqrt t$ (see, e.g., [KrS]). 
This is related to the conjectured
identity $\xi = 2\chi$ (for $d=2$) for the fluctuation exponents
$\xi$ and $\chi$ that are the main topic of the next section.
We remark that the exact values $\chi=1/3$ and $\xi=2/3$ have been
derived recently in [BDJ, J] for a model related to random
permutations, one of whose many guises
is a kind of $d=2$ directed FPP.

There are many interesting open questions one can ask about height
functions on $Q$ satisfying (1.11) and (1.12), such as whether there
exist ones essentially different from those of the form 
$H^{\hat x}(q,q(0))$.  For example, in general $d$ one could 
take two deterministic sequences of points $\bar q^{(1)}_k$ and $\bar q^{(2)}_k$
with $|\bar q^{(1)}_k| = |\bar q^{(2)}_k|\to\infty$
and with $\bar q^{(j)}_k/|\bar q^{(j)}_k|\to \hat x^{(j)}$ for
$j=1,2$ as $k\to\infty$ and then study
$$
\min_{j=1,2} (T_\a(q,\bar q_k^{(j)})) - 
\min_{j=1,2} (T_\a(q',\bar q_k^{(j)}))\tag1.15
$$
as $k\to\infty$.  It could be (and this seems likely the case
for $d=2$) that the limit (in distribution) of this random function
of $q$ and $q'$ is a symmetric mixture of the distributions of
$H^{\hat x^{(1)}}$ and $H^{\hat x^{(2)}}$.
This would be because the boundary between the region
of $Q$ where $T_\a(q,\bar q^{(1)}_k) < T_\a(q,\bar q^{(2)}_k)$
and where $T_\a(q,\bar q^{(1)}_k) > T_\a(q,\bar q^{(2)}_k)$
would (probably) be far from the origin as $k\to\infty$.
On the other hand, it is conceivable (e.g., for large
enough $d$, if $\chi=0$ there---see the discussion and references
in [KrS] or [NewP]) 
that this
boundary would not wander off to infinity but rather would
have an a.s. limit, and thus that 
(1.15) would also have a limit. 
The latter limit, defined for all $q,q'$, should
equal either $H^{\hat x^{(1)}}$ or $H^{\hat x^{(2)}}$, but only
when $q$ and $q'$ are both on the same side of the
limit boundary.
Thanks to Theorem 1.12, we can now pose such questions, but answering them
remains a task for the future.

\head{2.  Fluctuation Results}\endhead

Throughout this section and the remainder of the paper we deal
with some fixed $d\ge2$ and $\a\in(1,\infty)$.  Occasionally, as noted,
we will restrict our attention to $d=2$ and $\a\in[2,\infty)$.
We drop the $\a$ subscript in the (pseudo-) metric $T_\a(x,y) = 
T_\a(q(x),q(y))$ and the geodesic $\M(x,y) = \M(q(x),q(y))$
and denote these by $T(x,y)$ and $M(x,y)$.  This section is organized
as follows.  In Subsection 2.1, we state two
theorems  giving large
deviation bounds on $T(x,y)$ as $|x-y|\to\infty$; the proofs are given
later in Sections 3 and 4.  The first theorem
concerns fluctuations about the mean and the second concerns fluctuations
about $\mu|x-y|$.  Here $\mu = \lim_{|x-y|\to\infty} ET(x,y)/|x-y|$
and also equals the a.s. limit of $T(0,n\hat e)/n$ as $n \to \infty$
for any fixed unit vector $\hat e$ [HoN1]; it is of course the
same $\mu$ appearing in the Shape Theorem 1.4.
A third theorem in Subsection 2.1 gives a strengthened shape
theorem like the one obtained for lattice FPP in [Al3, Ke2].
In Subsection 2.2, we state and prove (using the theorem about
$T(x,y)-\mu|x-y|$) results about fluctuations of $M(x,y)$ from
a straight line as
$|x-y|\to\infty$.  These results tell us that, with high probability,
long finite geodesics (a) do not deviate far from the straight line between
their endpoints and (b) do not start off in one direction and
then ``noticeably'' change course.  In Subsection 2.3, we apply
these fluctuation results to prove Theorems 1.8-1.11.

\subhead{2.1 Fluctuation of the Metric}\endsubhead
In this subsection we consider fluctuations of $T_\ell \equiv T(0,\x)$
where $\ell>0$ and $\hat e_1$ is the unit vector $(1,0,\dots,0)$.
As in the case of lattice models, one expects that the standard deviation
of $T_\ell$ grows like $\ell^\chi$ for some exponent
$\chi = \chi(d)$ that should not depend on $\a$.  For lattice
FPP on ${\Bbb Z}^1$ (with $\ell$ an integer), the analog of $T_\ell$ 
is the sum of $\ell$ \iid\ random variables ($\tau(j-1,j)$ with $1\le j\le\ell$)
so that $\chi(1) = 1/2$ (assuming $E[\tau(j-1,j)^2] <\infty$). For
Euclidean FPP on ${\Bbb R}^1$, again $\chi(1) = 1/2$ although the
argument, while standard, is not as trivial since $T_\ell$ then is essentially
$\sum_{i=1}^N U_i^\a$ where the $U_i$ are \iid\ exponential
random variables and $N$ is random such that
$\sum_{i=1}^N U_i$ is close to $\ell$. For
$d=2$, $\chi(2)$ is believed to equal $1/3$ 
(see [HuH, K, HuHF, KPZ]), but the only
models for which this (and much more) has been proved are certain
directed FPP-like models related to random permutations (see [BDJ, J]). For
lattice FPP with 
$d\ge2$, there have been rigorous bounds on $\chi(d)$ including
Kesten's result that $\chi(d)\le 1/2$ [Ke2].   This latter bound
has been strengthened by Kesten [Ke2] and Alexander [Al1, Al3] 
to give large
deviation upper bounds for the deviation of $T_\ell$ as $\ell\to\infty$
from its mean and from the asymptotic expression $g(\x)$
(or more generally for the deviation of
$T(q,q')$ for $q,q'\in \zd$ as $|q-q'|\to\infty$
from its mean and from $g(q-q')$), where
$$
g(v)\qequals \lim_{n\to\infty} \frac{E[T(0,nv)]}{n}. \tag2.1
$$

In the case of lattice FPP, $g$ is a norm on $\rd$ whose
unit ball arises in the shape theorem [R, CD, Ke1, Bo].  For
Euclidean FPP,
$$
\lim_{\ell\to\infty} \frac{E[T(0,\ell\hat x)]}{\ell}
\qequals \lim_{\ell\to\infty}\frac{ET_\ell}{\ell} \qequals \mu
\ \in\ (0,\infty)\tag2.2
$$
(see (7) of [HoN1]),
where $\mu = \mu(\a,d)$  appears in the shape
theorem (Theorem 1.4) above.  
The next two theorems are the analogs for Euclidean
FPP of the Kesten and Alexander results of [Al3, Ke2] for lattice
FPP.  The great advantage of Euclidean FPP over the lattice case is
that the unit ball of the metric $g(v)$ (about which very little has
been proved) is replaced by the {\it Euclidean} ball (of radius $\mu^{-1}$).
This allows us in the next subsection to go well beyond what was proved for
lattice FPP, as we discuss there.

Here and for the remainder of the paper, we use $\cs$ 
to represent a strictly positive constant, to be thought of as small,
that depends on $\a$ and $d$ but never on $\ell$.  The actual value of $\cs$
may decrease as the paper progresses
(perhaps even in a single line);  all statements made involving $\cs$ are
valid with any smaller choice 
of $\cs$.  Analogously, $\cb$ is a positive finite constant,
thought of as large, whose value does not depend on $\ell$ but increases (with
similar impunity) as the
paper progresses.  
Certain other constants, appearing as exponents, we keep track
of more carefully.  We record their values here for easy reference:
$$
\aligned
\kf &\qequals \min(1,d/\a),\\
\kb &\qequals 1/(4\a+3),\\ 
\kc &\qequals 1/(2\a), \\
\kd &\qequals d/\a, \text{\ and} \\
\ke &\qequals 1/(4\a+2).\\
\endaligned\tag2.3
$$

\proclaim{Theorem 2.1}  Let $d\ge2$ and $\a>1$. For
some constant $\cb$, 
$$
\var T_\ell \qle \cb \ell\text{\ \, for}\ \,  \ell\ge0.\tag2.4
$$
Additionally, with $\kf=\min(1,d/\a)$, $\kb=1/(4\a+3)$, and 
for some constants $\cs$ and $\cb$, 
$$
\prob[\,|T_\ell-ET_\ell| > x\sqrt{\ell}] \qle \cb\exp(-\cs x^\kf)
\text{\ \, for}\ \,  \ell\ge0 \text{\ \, and}\ \,
0\le x\le \cs\ell^{{\kb}}.\tag2.5
$$
\endproclaim
The proof of Theorem 2.1 is given in Section 3.  The next theorem,
which is essentially a replacement of $ET_\ell$ by $\mu\ell$ in (2.5),
is proved in Section 4, by using Theorem 2.1 to show that
$$
|ET_\ell - \mu\ell| \qle \cb \sqrt{\ell} (\log\ell)^{1/\kf}.\tag2.6
$$

\proclaim{Theorem 2.2} Let 
$d\ge2$, $\a>1$, $\kf=\min(1,d/\a)$ and $\kb=1/(4\a+3)$. For 
any $\e$ in $(0,\kb)$, there exist constants $\cs$ and $\cb$ (depending on
$\e$) such that
$$
P[|T_\ell - \mu\ell| \ge\lambda] \qle \cb\exp(-\cs(\lambda/\sqrt{\ell})^\kf)
\text{\ \,for}\ \, \ell >0 \text{\ \,and}\ \, 
\ell^{\frac1 2 + \e} \le \l\le \ell^{\frac1 2 + \kb -\e}.
\tag2.7
$$
\endproclaim

A corollary of Theorem 2.2, the proof of which we
sketch in Section 4, is the following improvement of Theorem 1.4;
it is an analog of the Alexander-Kesten improved shape theorem
for lattice FPP [Al3]:

\proclaim{Theorem 2.3} For any $\a\in(1,\infty)$ and $d\ge 2$,
with ${\Cal B}_0 \equiv {\Cal B}(0,\mu^{-1})$, 
the following is true almost surely:
$$
Q\cap\Big(1 - \frac{(\log s)^{2/\kf}}{\sqrt{s}}\Big)s{\Cal B}_0 \ \subset\  B_\a(0,s^{1/\a})\ \subset
\ \Big(1 + \frac{(\log s)^{2/\kf}}{\sqrt{s}}\Big)s{\Cal B}_0
$$
for all sufficiently large $s$.
\endproclaim

We make no claims about the optimality of the exponents $\kf$ and
$\kb$ appearing in (2.5)-(2.7). We also note that the power $2/\kf$
in Theorem 2.3 can be replaced by $(1+\e)/\kf$ with any $\e > 0$. For
lattice FPP with an exponential tail
assumption on the underlying $\tau(q,q')$ variables, the analogous
results in [Al3,Ke2] have $\kf=1=\kb$. 
In the next subsection, we use Theorem 2.2 to control deviations
of long finite geodesics from approximately straight line behavior.

\subhead{2.2 Fluctuations of Geodesics}\endsubhead
We want to use Theorem 2.2 to
bound the probability that the geodesic $M(x,y)$ touches a 
Poisson particle located far from the straight line segment
$\overline{xy}$.  Our reasoning will follow that used in
[New1] (see (3.2) there) but modified for the
Euclidean context. 
We use (2.7) and some other arguments to show that for any $\e>0$,
with high probability for large $|x-y|$, $M(x,y)$ does not deviate
more than order $|x-y|^{\frac3 4 +\e}$ from $\overline{xy}$. 
The wandering exponent $\xi = \xi(d)$ may be regarded as defined
so that $|y-x|^\xi$ is the actual order of the typical 
(or largest) deviation
from $\overline{xy}$.  Thus, our next theorem implies that 
$\xi\le 3/4$.  It is conjectured that $\xi(2)=2/3$ and decreases
to $1/2$ for increasing $d$ (see the discussion and references in
[KrS] or [NewP]).  
A related result was obtained
in [NewP] for lattice FPP but it was much weaker because
of lack of information about the asymptotic shape ${\Cal B}_0$ for lattice
FPP.  Roughly speaking, the lattice result was only valid
when $y-x$ points in a direction where the boundary of ${\Cal B}_0$
is curved.  If it were proved 
that in a lattice model ${\Cal B}_0$ is uniformly curved, then a
result like the next theorem (which is only for {\it Euclidean\/} FPP) would
follow --- see [NewP] for details.

We define $M_\ell = M(0,\x)$ and, for $A\subset\rd$,
$$
\dmax(M_\ell,A) \qequals \sup_{q\in M_\ell} \dist(q,A),
\tag2.8
$$
where $\dist(q,A)$ denotes the ordinary Euclidean distance from $q$ to the set $A$.
This represents the maximal Euclidean distance of 
(any point in) $M_\ell$ from $A$;
if $M_\ell$ is replaced by a single point $y$, then $d_{\text{max}}(y,A)$
is the usual Euclidean distance of $y$ to the set $A$.

\proclaim{Theorem 2.4}  Let $d\ge2$, $\a>1$, $\kf=\min(1,d/\a)$ 
and $\kb=1/(4\a+3)$. For any $\e\in (0,\kb/2)$, there exist
constants $\cs$ and $\cb$ (depending on $\e$) such that
$$
P[d_{\text{max}}(M_\ell,\overline{0\,\x})\ge\ell^{\frac3 4 + \e}]
\qle \cb\exp(-\cs\ell^{\,3\e\kf/4}).
\tag2.9
$$
Furthermore, with $\B = \B(\x,1) = \{x\in\rd:|x-\x|\le1\}$, 
for (possibly different) $\cs$ and $\cb$,
$$
P[\exists\ b\in \B\text{\ with }
\dmax(M(0,b),\overline{0\,b})\ge|b|^{\frac3 4 + \e}]
\qle \cb\exp(-\cs\ell^{\,3\e\kf/4}).
\tag2.10
$$

\endproclaim
\demo{Proof} We will prove that, for some $\cs$ and $\cb$,
$$
P[\exists\ b\in \B\text{\ with }\dmax(M(0,b),\overline{0\,\x})\ge\ell^{\frac3 4 + \e}]
\qle \cb\exp(-\cs\ell^{\,3\e\kf/4}),\tag2.11
$$
from which (2.9) follows immediately and (2.10) follows (for possibly different $\cs$ and $\cb$)
from the facts that $|\,|b|-\ell| \le 1$
and $|\dmax(M(0,b),\overline{0\,b}) - \dmax(M(0,b),\overline{0\,\x})| \le 1$.

We begin with the observations that 
$$
|T(u,v) -T(u,w)| \qle |q(v) - q(w)|^\a \quad \text{for all $u,v,w\in\rd$}, \tag2.12
$$
and that, for $q\in Q$ and $w\in\rd$, $|q-q(w)|\le 2|q-w|$, so also
$$
|T(u,q) - T(u,w)| \qle (2|q-w|)^\a\quad\text{for all $q\in Q$ and $u,w\in\rd$.}\tag2.13
$$
Furthermore, repeated application 
of the triangle inequality to (2.12) gives that
$$
|T(u,v) -T(u,w)| \qle (2|q(v)-v|+2|v-w|)^\a \quad \text{for all $u,v,w\in\rd$}. \tag2.14
$$

Now let
$$
\align 
A'_\ell &\qequals \{x\in\rd : \dist(x,\overline{0\,\x}) < \ell^{\frac3 4 +\e}\},\\
A_\ell  &\qequals \{x\in \rd\setminus A'_\ell : \dist(x,A'_\ell) < \ell^{\frac3 4}\},\text{\ and}\\
A^+_\ell  &\qequals \{x\in \rd\setminus A'_\ell : \dist(x,A'_\ell) < \ell^{\frac3 4}+\sqrt{d}\}.
\endalign
$$
Additionally, let $F_\ell$ denote the event
that $q(0), q(\x)\in A'_\ell$ and every geodesic
segment $\overline{q_kq_{k+1}}$ with either $|q_k|\le\ell$ or
$|q_{k+1}|\le\ell$ has $|q_k - q_{k+1}|\le\ell^{3/4}$.
By an application of Lemma 5.2 (see (5.5)),   $F_\ell$
satisfies $P[F_\ell^c]\le \cb\exp(-\cs\ell^{3d/4})$.
Furthermore, for large $\ell$, on $F_\ell$, for $b\in \B$ we have
$$
\align 
\dmax(M(0,b),\overline{0\,\x}) \ge \ell^{\frac3 4+\e} 
\ \implies\ &\text{$\exists$ $q\in Q\cap A_\ell$ on $M(0,b)$} \\
\ \implies\ &\text{$\exists$ $q\in Q\cap A_\ell$
with $T(0,q) + T(q,b) = T(0,b)$}\\
\ \implies\ &\text{$\exists$ $w\in \zd\cap A^+_\ell$ with }
T(0,w) + T(w,\x) \le\\ &T(0,\x) +2((2\sqrt{d})^\a + (2|q(\x)-\x|+2)^\a).\\
\endalign
$$
This latter implication uses (2.13) and (2.14).
It follows that, on $F_\ell\cap\{|q(\x)-\x| < \ell^{1/(2\a)}\}$, 
for large $\ell$ and $b\in \B$,
$$
\align
\dmax(M(0,b),\overline{0\,\x}) \ge \ell^{\frac3 4+\e} 
\ \implies\ &\text{$\exists$ $w\in \zd\cap A^+_\ell$ with}\\
&T(0,w) + T(w,\x) \le T(0,\x) +\ell^{\frac1 2 + \e}.\\
\endalign
$$
Hence
$$
\aligned
P[\exists\ b\in \B\text{\ with }&\dmax(M(0,b),\overline{0\,\x})\ge\ell^{\frac3 4 + \e}]\\
\qle&\cb\exp(-\cs\ell^{3d/4})\ +\ \cb\exp(-\cs\ell^{d/(2\a)})\\
&+\ \sum_{w\in\zd\cap A^+_\ell} P[T(0,w) + T(w,\x) \le T(0,\x) +\ell^{\frac1 2 + \e}].\\
\endaligned\tag2.15
$$
The proof is completed by combining (2.7) and
(2.15) with some elementary geometry.

Given $x,y\in\rd$, let $\Delta(x,y) = \mu(|y| + |x-y| - |x|)$.  Then on
$\{T(0,w) + T(w,\x) \le T(0,\x) + \ell^{\frac1 2 + \e}\}$ at least one of
$\big|T(0,w) - \mu|w|\big|$, $\big|T(w,\x) - \mu|\x-w|\big|$, or $\big|T(0,\x) - \mu\ell\big|$
must exceed $\tilde\Delta(\x,w)\equiv(\Delta(\x,w) - \ell^{\frac1 2 + \e})/3$, so
$$
\aligned
P[T(0,w) + T(w,\x) \le T(0,\x) +\ell^{\frac1 2 + \e}]\qle &P[\big|T(0,w) - \mu|w|\big|>\tilde\Delta(\x,w)]\\
 &+\ P[\big|T(w,\x) - \mu|\x-w|\big| > \tilde \Delta(\x,w)]\\
 &+\ P[\big|T(0,\x) - \mu\ell\big| >\tilde\Delta(\x,w)]\\
\endaligned\tag2.16
$$

We note that, for $w\in A_\ell^+$, $\Delta(\x,w)$, and hence $\tilde\Delta(\x,w)$, 
is at least of order
$\ell^{\frac1 2 + 2\e}$ and at most of order $\ell^{\frac3 4 + \e}$
as $\ell\to\infty$. For example, for $w = \x/2 + (\ell^{\frac 3 4 + \e}+\sqrt{d})\hat e_2$,
$\Delta(\x,w)/\mu = 2((\ell/2)^2+(\ell^{\frac3 4 +\e}+\sqrt{d})^2)^{1/2} -\ell$ 
which is of order
$\ell^{\frac1 2 + 2\e}$ by the Pythagorean theorem,
while, for $w = (-\ell^{\frac 3 4 + \e}-\sqrt{d})\hat e_1$, 
$\Delta(\x,w) = 2\mu(\ell^{\frac3 4 + \e}+\sqrt{d}) = O(\ell^{\frac3 4 + \e})$.  

Each of the three probabilities in (2.16) 
can be expressed (using Euclidean invariance)
in the form of the probability of (2.7) with $\ell$
replaced by some $\ell'$ between order $\ell^{\frac 3 4 + \e}$ and
order $\ell$, and with $\lambda$ between order
$\ell^{\frac1 2 + 2\e}$ and order $\ell^{\frac 3 4 + \e}$.
We can bound these probabilities by replacing $\l$ by the smaller
$\l' = (\ell')^{\frac1 2 + \e}$ so that the
condition $(\ell')^{\frac1 2 + \e} \le \lambda' 
\le (\ell')^{\frac1 2 + \kb - \e}$ is satisfied.  Since 
$A_\ell$ can be contained in a Euclidean ball of radius
order $\ell$, we have 
$$
\align
P[\exists&\ b\in \B\text{\ with }\dmax(M(0,b),\overline{0\,\x})\ge\ell^{\frac3 4 + \e}]\\&\qle
\cb\exp(-\cs\ell^{3d/4}) 
+\ \cb\exp(-\cs\ell^{d/(2\a)})
+\ \cb\ell^d\sup\{\exp(-\cs((\ell')^{\frac1 2 + \e}/\sqrt{\ell'})^{\kf})\}, \tag2.17\\
\endalign 
$$
where the supremum is over all $\ell'$ with
$\ell^{\frac3 4+\e}\le\ell'\le\ell$.
This yields (2.11) by taking 
$\ell' = \ell^{\frac3 4+\e}$ and noting that for large $\ell$,
$(\ell')^{\e \kf} \ge \ell^{(\frac3 4 + \e) \e \kf}
\ge \ell^{3\e\kf/4}$.\qed
\enddemo

Theorem 2.4 immediately yields:

\proclaim{Corollary 2.5} For $d\ge2$, $\a>1$ 
and any $\e>0$, the number $N_\e$ of
$q'\in Q$ such that 
$\dmax(M(0,q'), \overline{0\,q'}) \ge |q'|^{\frac3 4 + \e}$
is \as\ finite.
\endproclaim

\demo{Proof} It follows from (2.10) of Theorem 2.4, rotational invariance,
and an application of the Borel-Cantelli
Lemma, that \as\ the events  
$$
\{\exists\ b\in \B(w,1)\text{\ with }\dmax(M(0,b),\overline{0\,b})\ge|b|^{\frac3 4 + \e}\}
$$
occur for only finitely many $w\in(2/\sqrt{d})\zd$.  
The corollary follows since the $\B(w,1)$ cover $\rd$. \qed
\enddemo

The next theorem, which itself is a consequence of this corollary, gives a different version of
the inequality $\xi\le 3/4$.  To formulate the theorem, we need some notation.  Let
$C(x,\e)$ for nonzero $x\in\rd$ and $\e\in[0,\pi)$ denote the cone
$$
C(x,\e) \ \equiv\ \{y\in\rd:\theta(x,y) \le\e\},\tag2.18
$$
where $\theta$ is the angle (in $[0,\pi]$) 
between $x$ and $y$.
Recalling the definition of the spanning tree $R(q) = R_\a(q)$ formed by
all geodesics $M(q,q')$ from $q$ as given in Proposition 1.2, we define
$R^{\text{out}}(q,q')$ for $q'\in Q$ to be the set of all $q''\in Q$ such that
$M(q,q'')$ touches $q'$ --- i.e., it is the
part of $R(q)$ going ``outward'' from
$q'$.  The next theorem states that for any $q$ and all but 
finitely many $q'$ (the number depending on $q$),
any geodesic continuation of $M(q,q')$ must remain 
inside $q + C(q'-q,\h^*(|q'-q|))$ where 
$\h^*(\ell) \equiv \ell^{\frac3 4 + \e} / \ell$.  
This was announced as Theorem 2 of [HoN2].

\proclaim{Theorem 2.6} Let 
$d\ge2$, $\a>1$, $\e\in(0,\frac1 4)$, and $\h^*(\ell) = \ell^{-\frac1 4 + \e}$.
Then almost surely, for every $q\in Q$, for all 
but finitely many $q'\in Q$,
$$
R^{\text{out}}(q,q') \ \subset\ q + C(q'-q, \h^*(|q'-q|)). \tag2.19
$$
\endproclaim

\demo{Proof} It suffices to restrict attention to $q=q(0)$. From 
Lemma 5.2  (see (5.5)) and the Borel-Cantelli
Lemma, it follows that there is some finite (random) $L_0$ so 
that for any geodesic segment
$\overline{q_k\,q_{k+1}}$ 
with $|q_k|\ge L_0$, $|q_{k+1}-q_k| \le |q_k|^{3/4}$. Theorem 2.6
is then a consequence of Corollary 2.5 and the following deterministic lemma.
\enddemo

\proclaim{Lemma 2.7} Let $d\ge2$ and $\d\in(0,\frac1 4 )$. Suppose
$(q_i) = (q_1,q_2,\dots)$ is any sequence of distinct points
in $\rd$ with $|q_i|\to\infty$ such that for all large $j$,
$$
|q_{j+1}-q_j|\le|q_j|^{3/4}\text{\ \ and\ \ }
\dist(q_k,\,\overline{0\,q_j})\le |q_j|^{1-\d}
\text{\ \ for $1\le k< j$.}\tag2.20
$$
Then there exists $\cb$ and $k^*>0$ such that
$$
\theta(q_k,q_j) \qle \cb |q_k|^{-\d}
\text{\ \ whenever $k^*\le k<j$}.\tag2.21
$$
\endproclaim

\demo{Proof} Choose $L$ large enough that $L^{3/4} < L^{1-\d} < L/3$,
and then choose $k^*$
so that (2.20) holds and $|q_j|\ge L$ whenever $j\ge k^*$.
Now suppose $k^*\le k<j$.

{\it Case 1:  $|q_j| \le 3|q_k|$.}  First note that we must have
$\theta(q_k,q_j) < \pi/2$, for otherwise
$$
\dist(q_k,\overline{0\,q_j}) \qequals |q_k|
\qge \frac{|q_j|}{3} \ >\ |q_j|^{1-\d},
$$
which violates the second part of (2.20). It follows then from 
elementary geometric considerations that
$$
\sin \theta(q_k,q_j) \qle \frac{|q_j|^{1-\d}}{|q_k|}
\qle 3^{1-\d}|q_k|^{-\d}.
$$
Using that $\theta\le\frac{\pi}{2}\sin\theta$ on $[0,\frac{\pi}{2})$,
we see that $\theta(q_k,q_j)\le\cb|q_k|^{-\d}$.

{\it Case 2: $|q_j| > 3|q_k|$.} We will construct a subsequence
$(q_{i_0},\dots,q_{i_n})$ of $(q_k,\dots q_j)$ such that:
$q_{i_0} = q_k$; $q_{i_n}=q_j$; 
the $|q_{i_m}|$ are increasing; $|q_{i_{m+1}}|\le 3|q_{i_{m}}|$
for $m+1 \le n$;
and, for $m+1 \le n-1$, $|q_{i_{m+1}}| \ge 2|q_{i_m}|$.
As we shall presently see, this is possible
because, by the first part of (2.20), 
the sequence $(q_i)$ stretches from $q_k$ to $q_j$ without
any (relatively) large gaps.
We then will have 
$|q_{i_m}| \ge 2^{m-1}|q_{i_0}|=2^{m-1}|q_{k}|$ for $0\le m\le n$,
with the exponent $m-1$ (instead of $m$)
to accomodate the case $m=n$. It
follows from this and a repeated application of Case 1 that
$$
\aligned
\theta(q_k, q_j)\qequals\theta(q_{i_0},q_{i_n}) 
&\qle \sum_{m = 0}^{n-1}\theta(q_{i_m},q_{i_{m+1}})
 \qle \sum_{m=0}^{n-1} \cb |q_{i_m}|^{-\d}\\
&\qle \cb\bigg(\sum_{m = 0}^{n-1} 2^{-(m-1)\d} \bigg)|q_k|^{-\d}
\qle \cb |q_k|^{-\d},\\
\endaligned
$$
where the final inequality holds for a larger $\cb$.

To construct the requisite $(q_{i_m})$, put
$i_0 = k$ and suppose $i_m$ has been selected.
If $i_m = j$, put $n=m$ and stop.  Otherwise, let
$i_{m+1} = \max\{i: i_m < i \le j,\ |q_i| \le 3|q_{i_m}|\}$.
By construction, the $|q_{i_m}|$ are increasing with 
$|q_{i_{m+1}}|\le 3|q_{i_m}|$.  Furthermore,
for $m+1 \le n-1$ (so also $i_{m+1} < j$), 
we must have $2|q_{i_m}|\le|q_{i_{m+1}}|$, for otherwise
$$
\aligned
|q_{i_{m+1}+1} - q_{i_{m+1}}|&\qge |q_{i_{m+1}+1}| - |q_{i_{m+1}}|\\
&\qg  3|q_{i_m}| - 2|q_{i_m}| \qequals |q_{i_m}| 
\qge \frac {|q_{i_{m+1}}|}{3} \qge |q_{i_{m+1}}|^{3/4},\\
\endaligned
$$
in contradiction to the first part of (2.20).\qed

\enddemo

\subhead{2.3 Proof of Theorems 1.8-1.11}\endsubhead
Suppose $R$ is a tree
whose vertex set is an infinite countable subset of $\rd$ with
$u$ and $u'$ two vertices of $R$.  We define $R^{\text{out}}(u,u')$,
as in the last subsection, 
to be the set of vertices $u''$ of $R$ such that the (unique) path
in $R$ from $u$ to $u''$ touches $u'$. 

\proclaim{Definition}  For $\h$ a positive function on
$(0,\infty)$, we say that such a tree 
$R$ is $\h$-straight at the vertex $u$ if for all but
finitely many vertices $u'$ of $R$,
$$
R^{\text{out}}(u,u')\ \subset\ u  + C(u'-u,\h(|u'-u|)).\tag2.22
$$
\endproclaim

Theorem 2.6 is the statement that \as, for every $q\in Q$,
$R(q)$ is $\h^*$-straight for $\h^*(\ell)=\ell^{-\frac1 4 + \e}$.

\proclaim{Definition}  $Q'$, a subset of $\rd$, is said to be
asymptotically omnidirectional if for
all finite $K$, 
$\{q/|q| : q\in Q' \text{\ and\ } |q|>K\}$ is
dense in $S^{d-1}$.
\endproclaim

\proclaim{Proposition 2.8} Suppose $R$ is a 
tree whose vertex set $U\subset\rd$ is
locally finite but asymptotically omnidirectional
and such that every vertex has finite degree.  Suppose further that
for some $u\in U$, $R$ is $\h$-straight at $u$, where $\h(\ell)\to0$ as
$\ell\to\infty$. Then $R$ satisfies the following properties:
(i) every semi-infinite path in $R$ starting from $u$ has an\
asymptotic direction; (ii) for every $\hat x\in S^{d-1}$, there
is at least one semi-infinite path in $R$ starting from $u$ with
asymptotic direction $\hat x$; (iii) the set $V(u)$ of
$\hat x$'s such that there is more than one semi-infinite
path starting from $u$ with asymptotic direction $\hat x$ is
dense in $S^{d-1}$.
\endproclaim

\demo{Proof} Let $u=u_1,u_2,\dots$ be a semi-infinite path
in $R$.  Then $\h$-straightness implies that for large $m$, the
angle $\theta(u_n-u,u_m-u) \le \h(|u_m-u|)$ for $n\ge m$. 
Since $|u_m|\to\infty$ as $m\to\infty$ (because $U$ is
locally finite), it follows that $u_n/|u_n|$ converges, proving
(i).  Since $U$ is asymptotically omnidirectional and
each vertex has finite degree, it follows
that starting from $v_1=u$, one can for a given $\hat x$ inductively
construct a semi-infinite path $v_1,v_2,\dots$ in $R$ such that
for each $j$, $R_{\text{out}}(u,v_j)$ contains a sequence (depending
on $j$) $u_1,u_2,\dots$ with $u_n/|u_n|\to\hat x$.  But (i) shows
that $v_j/|v_j|$ tends to some $\hat y$ and then $\h$-straightness
implies $\theta(\hat x,v_j-u)\le \h(|v_j-u|)$ so that letting
$j\to\infty$ yields $\hat x = \hat y$, proving (ii).

Given any (large) finite $K$, one can consider those (finitely many) vertices
$v$ with $|v|>K$ such that no other vertex $w$ on the path
from $u$ to $v$ has $|w|>K$. By taking a subset of these $v$'s,
one obtains a finite set
of vertices $v_1^{(K)},\dots,v_{m(K)}^{(K)}$ with
$|v_j^{(K)}|>K$ such that the $R^{\text{out}}(u,v_j^{(K)})$'s
are disjoint and their union includes all but finitely many
vertices of $U$.  For a given $K$, let $G_j$ denote the set of $\hat x$'s
such that some semi-infinite path from $u$  passing through
$v_j^{(K)}$ has asymptotic direction $\hat x$.
Then by (ii), $\cup_j G_j = S^{d-1}$.  On the other hand,
by $\h$-straightness, each $G_j$ is a subset of the
(small) spherical cap 
$\{\hat x:\theta(\hat x, v_j^{(K)}) \le \h(|v_j^{(K)}-u|)\le \e(K)\}$
where $\e(K)\to0$ as $K\to\infty$ (since $|v_j^{(K)}| > K)$.
Furthermore, by the same arguments that proved (ii), each $G_j$ 
is a {\it closed\/} subset of $S^{d-1}$.  It follows that $V(u)$
contains, for each $K$, $\cup_{j\le m(K)} \partial G_j$,
where $\partial G_j$ denotes the usual boundary ($G_j$ less its interior).
Since $\e(K)\to0$ as $K \to \infty$, we obtain (iii)
by standard arguments.\qed
\enddemo

\demo{Proof of Theorems 1.8, 1.9 and 1.10} These three theorems
are all essentially immediate
consequences of Proposition 2.8
and the (easily proven) fact that $Q$ is \as\ locally
finite and asymptotically omnidirectional.
\enddemo

\demo{Proof of Theorem 1.11} The only part of Theorem 1.11 that
remains to be proved (i.e., that does not immediately follow from
Theorems 1.6 and 1.8) is that $(\hat x,\hat y)$-geodesics
with $\hat y \not=-\hat x$ do not exist, even for nondeterministic
$\hat x$ and $\hat y$ depending on $Q$.  To prove this, it
suffices to show, for each $\d>0$, that this is the case with
the further restrictions that $\theta(\hat y, -\hat x) >\d$ and
that the $(\hat x,\hat y)$-geodesic touches $q(0)$.  Let 
$E_k$ denote the event that there exist $q,q''\in \rd$ with
$|q'|,|q''|\in[k,k+k^{3/4}]$, $\theta(q'',-q')>\d/2$, and with
$M(q',q'')$ touching $q(0)$.  By arguments like those in the
proofs of Theorem 2.3 and Corollary 2.4 one can prove that
$P[E_k \text{\ infinitely often}] = 0$ and that this leads
to the non-existence of $(\hat x,\hat y)$-geodesics
passing through $q(0)$ with $\theta(\hat y,-\hat x)>\d$.
We leave further details to the reader.\qed
\enddemo

\head{3.  Proof of Theorem 2.1}\endhead

In many respects, our proof of 
Theorem 2.1 parallels the arguments in [Ke2], where
analogous results for lattice FPP are presented.  However, our
Euclidean framework presents a host of technical issues. For
such technical reasons 
we will need to work with certain approximations of
$T_\ell$.  With $\bar Q$ any locally finite subset of $\rd$,
$\phi:{\Bbb R}^+\to{\Bbb R}^+$ any continuous
strictly increasing convex function with $\phi(0)=0$ (the cost function)
and $a$ and $b$ arbitrary, and possibly random, points in $\rd$, define
$$
T[\bar Q, \phi, a, b] \qequals \inf\bigg\{
\sum_{j=1}^{k-1} \phi(|q_j-q_{j+1}|): 
k\ge2, q_1 = a, q_k = b, q_j\in\bar Q \text{ for } 1<j<k\bigg\}.\tag{\eqab}
$$
So, for example, with $\phi_\infty(t) \equiv t^\alpha$, 
we have
$T_\ell = T[Q, \phi_\infty, q(0), q(\ell\hat e_1)]$.
Our first approximation to $T_\ell$, denoted by $T'_\ell$, is defined by
$$
T'_\ell \qequals T[Q, \phi_\infty, 0, \x].
$$
It would seem that $T'_\ell$ is a more natural quantity to study, 
since the paths under consideration
actually start at $0$ and end at $\x$.  Unfortunately, $T'_\ell$ does not obey a triangle 
inequality whereas $T_\ell$ does. For 
our second approximation, $T_\ell''$, 
we will need a collection
of subsets $Q_\ell\subset Q$ to be defined later (see above (\eqaf)) and a
family of cost functions $\phi_\ell$ 
defined as
$$
\phi_\ell(t) \qequals \cases 
                      t^\alpha &\text{if $t \le h_\ell$}\\
                      h_\ell^\alpha + \alpha h_\ell^{\alpha-1}(t-h_\ell)&\text{otherwise},
                      \endcases \tag\eqbu
$$
where $h_\ell= \max(h_0, h_1\ell^{{\kc}})$ with ${\kc} = 1/(2\a)$, 
and with both $h_0\ge1$ and $h_1\ge h_0$ to 
be specified later (see above (\eqaw) and below (\eqbnew)).  
Note that $\phi_\ell(t)\uparrow\phi_\infty(t) = t^\a$
as $\ell\to\infty$;  we will also have $Q_\ell\uparrow Q$.
We now define 
$$
T_\ell''\qequals T[Q_\ell, \phi_\ell, 0, \ell\hat e_1].
$$
These approximations to $Q$ and $\phi_\infty$ will play
a role similar to a truncation argument allowing $T''_\ell - ET''_\ell$ to be expressed
as the limit of a martingale with bounded differences.  

Throughout this section, we use the following notation. We let
$$
\aligned
q(0)&= r_1,r_2,\dots,r_K=q(\x),\cr
  0 &= r'_1,r'_2,\dots,r'_{K'}=\x,\text{\ and}\cr
  0 &= r''_1,r''_2,\dots,r''_{K''}=\x\cr
\endaligned
$$
achieve the infima in (\eqab) 
corresponding to $T_\ell$, $T'_\ell$, and $T''_\ell$ respectively.
We use $L_k$ to denote the ``link'' (i.e., straight line segment)
$\overline{r_kr_{k+1}}$, and we use $\overline r$ to denote the
polygonal path $L_1L_2\dots L_{K-1}$ with analogous interpretations of 
$L'_k$, $L''_k$, $\overline {r'}$, 
and $\overline {r''}$. For any link $L$, $|L|$ will be
its Euclidean length.
Also, for any cost function $\phi$ of the form (\eqbu) and $a,b\in\rd$, let
$$
\W_\phi(a,b) \qequals \{c\in\rd: \phi(|a-c|) + \phi(|c-b|) \le \phi(|a-b|)\}
$$
and put $\W(a,b) = \W_{\phi_\infty}(a,b)$.  
A number of properties of these subsets of $\rd$ are gathered in Lemma 5.1
of Section 5 below and used in the proof of the next lemma. 

With an appropriate $Q_\ell$ and $h_\ell$
the random variables $T_\ell$, $T_\ell'$, and $T_\ell''$ are related as follows:

\proclaim{Lemma 3.1} With ${\kd} = d/\a$ 
and for some constants $\cs$ and $\cb$, 
$$
\align
P[|T_\ell-T_\ell'|>x] &\qle \cb\exp(-\cs x^{{\kd}}), \text{\ for}
\ x>0 \tag{\eqac}\\
\text{\ and}\\
P[T_\ell'\not=T_\ell''] &\qle \cb\exp(-\cs\ell^{{\kc}}).\tag{\eqad}\\
\endalign
$$
\endproclaim

\demo{Proof of (\eqac)} The left 
side of (\eqad) is ill-defined until the $Q_\ell$ and $h_\ell$
are chosen;  we defer its proof.  
This does not apply to inequality (\eqac), which is easier
to prove.  
Let $\G(a) = \sup\{|c-a| : \W(a,c)\cap Q = \emptyset\}$,
and set $\G_\ell = \G(0) + \G(\x)$.
Then
$$
T'_\ell \qle T_\ell + |q(0)|^\a + |q(\x) - \x|^\a \qle T_\ell + \G_\ell^{\,\a},
$$
and, similarly, on $\{K'\ge3\}$,
$$
T_\ell\qle T'_\ell + |q(0)-r'_2|^\a + |q(\x)-r'_{K'-1}|^\a
 \qle T'_\ell + 2^\a \G_\ell^{\,\a},
$$
while on $\{K'=2\}$, $\G(0)\ge \ell$ so
$$
T_\ell \qle (|q(0)| + \ell + |q(\x)-\x|)^\a
\qle (\G(0)+\G(0)+\G(\x))^\a 
\qle 2^\a \G_\ell^{\,\a} \qle T'_\ell + 2^\a \G_\ell^{\,\a}.
$$
Collectively these bounds yields $|T_\ell - T'_\ell| \le 2^\a \G_\ell^{\,\a}$.
We complete the proof of (\eqac) by using the Remark following Lemma 5.2 below
(see (5.4)) to conclude that, for appropriate $\cs$ and $\cb$, 
$$
P[2^\a \G_\ell^{\,\a}>x] \qle P[\G(0) > x^{1/\a}/4]
+ P[\G(\x) > x^{1/\a}/4] \qle \cb\exp(-\cs x^{d/\a}).\tag{\eqae}
$$
\enddemo 

Our proof of (2.4) divides into the two cases $0\le \ell\le 1$ and $\ell>1$.
We are really only interested in the second (much more difficult) case, but proving the
first case illustrates the sort of technical difficulties created by our
definition of $T_\ell$.  We have: 
\proclaim{Lemma 3.2}  For some constant 
$\cb$, $\var T_\ell \le \cb\ell$ whenever $\ell\le1$.
\endproclaim
\demo{Proof}  If we were working 
with $T_\ell'$ instead of $T_\ell$, this case would be straightforward
since, for $\ell\le1$, $(T'_\ell)^2\le \ell^{2\alpha}\le  \ell$.  On the other
hand, although $T_\ell = 0$ 
for $\ell$ small enough that $q(0)=q(\x)$, no matter how small
$\ell$ is, among those Poisson particle configurations
for which $q(0)\not=q(\x)$, $|q(0)-q(\x)|$ (and $T_\ell$)
can be arbitrarily large. 
Looking a little closer, for any fixed $\ell\le1$ let $\tD$ denote the event
$\{q(0)\not=q(\x)\}$.  For $\rho>0$, on $\{|q(0)|=\rho\}$ we have
$$
T_\ell^2 \qle |q(0)-q(\x)|^{2\alpha} I_{\tD} \qle (|q(0)| + \ell
+ |q(\x)|)^{2\alpha} I_{\tD} \qle (2\rho+2)^{2\alpha}I_{\tD},
$$
where $I_{\tD}$ denotes the indicator of the event $\tD$.
Letting $A_\rho$ denote the event that there is a particle in the annulus
$\{x\in \rd:\rho < |x| < \rho+2\ell\}$,  we have
$$
\{|q(0)|=\rho\} \cap \tD \ \subset \ \{|q(0)|=\rho\} \cap A_\rho,
$$
so
$$
E[T_\ell^2\big||q(0)|=\rho] \qle (2\rho+2)^{2\a} P[A_\rho\big||q(0)|=\rho]\qle \cb (2\rho+2)^{2\a}(\rho+2)^{d-1}\ell,
$$
and 
$$
\align
\var T_\ell \qle ET_\ell^2 &\qequals \int_{\rho\ge0} E[T_\ell^2\big||q(0)|=\rho] dP[|q(0)|\le \rho]\\
& \qle   \ell \cb\int_{\rho\ge0} (2\rho+2)^{2\a+d-1} dP[|Q(0)|\le \rho]
\qequals \ell \cb.
\endalign
$$
(Recall that according to our conventions, the two instances of $\cb$
in the preceding equation represent different constants.) \qed

Proceeding with the case $\ell>1$, we define:
$$
S_\ell''\qequals \sum_{j=1}^{K''-1} \phi_\ell^2(|L''_j|).
$$
We do the case $\ell>1$ in three steps.

{\sl Step 1:  For any $\ell>0$,\ } $\var T''_\ell \qle 2^{2\a+1} \ex S_\ell''$.
We note that this inequality is also ill-defined until the $Q_\ell$ and  $h_\ell$ are
specified.  We presently define the $Q_\ell$; it turns out that Step 1 holds for any 
$h_\ell$.
Throughout this paper, for any length 
$\eta>0$, the ``$\eta$-boxes'' will refer to
the interior-disjoint
$d$-dimensional cubes whose vertices collectively 
are $\eta\cdot(\zd+(\frac1 2,\dots,\frac1 2))$. For
any $\eta$, the $\eta$-boxes 
may be associated with the $\zd$ lattice in the natural way:
$\nu\in\zd$ is associated with the box centered at $\eta\nu$.
Two $\eta$-boxes are called adjacent if they share a common $(d-1)$-dimensional face (i.e., if
their associated sites in $\zd$ are nearest neighbors). For
any Borel subset $S\subset\rd$,
let $\field(S)$ denote the $\s$-subfield of $\field$ generated by all events of the form
$\{\omega:Q(\omega)\cap B \not= \emptyset\}$ where $B$ ranges over all Borel subsets of $S$.
Clearly we may (and do) replace $\field$ 
with the possibly smaller $\field(\rd)$.  
Now fix any $\ell>0$ and
let $(B_m : m = 1 ,2, \dots)$ denote the $(\e/3^{\lfloor\ell\rfloor})$-boxes ($\e$ is
a quantity that depends only on $d$ and is specified in Step 2 below) enumerated in some order.  
We note that, in general, if $\eta'$ is an odd integral 
multiple of $\eta$ then the
$\eta$-boxes are nested in the $\eta'$-boxes so, in particular, the $(\e/3^{\lfloor\ell\rfloor})$-boxes 
are nested in the $\e$-boxes.
Let $q_m$ denote the leftmost particle in $Q\cap B_m$ (provided such a
particle exists) and let $Q_\ell = \{q_m\}\subset Q$ 
denote the set of all such leftmost particles.

Let $\field_m = \field(B_1\cup\dots\cup B_m)$ with
$\field_0 = \{\emptyset, \Omega\}$, so $\field_m\uparrow\field$
as $m\to\infty$.  Set
$$
\Delta_m \qequals \ex(T''_\ell\vert \field_m) - \ex(T''_\ell\vert \field_{m-1})
$$
so that
$$
T''_\ell - \ex T''_\ell \qequals \sum_{m = 1}^\infty \Delta_m, \text{\ and }
\var T''_\ell \qequals \sum_{m=1}^\infty \ex\Delta_m^2.
$$
This holds since $T''_\ell$ is bounded by $\ell^\a$.
Now set $\tilde\field_m = \field(\rd\setminus B_m)$ and
define
$$
\tilde\Delta_m \qequals T''_\ell - \ex(T''_\ell\vert \tilde\field_m).
$$
Then we have that
$\ex\Delta_m^2 \qle \ex\tilde\Delta_m^2$
since $\Delta_m = \ex(\tilde\Delta_m\vert\field_m)$ giving that
$$
\var T_\ell'' \qle \sum_{m=1}^\infty \ex\tilde\Delta_m^2
	\tag{\eqaf}
$$

In general, if $X$ and $Y$ are $L^2$ random variables with $Y$ measurable with
respect to some $\sigma$-field $\Cal G$, then 
$$
E[(X - E[X|{\Cal G}])^2 |{\Cal G} ] \qle E[(X-Y)^2|{\Cal G}].\tag\eqbx
$$
Put
$T_\ell^{(m)} = T[Q_\ell\setminus B_m, \phi_\ell, 0, \x]$; so $T_\ell^{(m)}$ 
is the minimal passage time from $0$ to $\x$
with respect to the $\phi_\ell$ cost function
using points in $Q_\ell$ other than in $B_m$, and
$T_\ell^{(m)}$ is  $\tilde\field_m$-measurable.  
Hence, with $U_m = (T_\ell^{(m)} - T''_\ell)^2$
we have
$$
E[\tilde\Delta_m^2|\tilde\field_m] \qle E[U_m|\tilde\field_m],
\text{\ and }
E\tilde\Delta_m^2 \qle EU_m. \tag\eqag
$$

Let $\bar{R}_m$ denote the event that  $q_m$
exists and equals $r''_i$ for some $i$.  
On the event $\bar{R}_m$, define the random variable $k(m)$ by the
relation $r''_{k(m)} = q_m$.  (Off of $\bar{R}_m$, the value of $k(m)$ is irrelevant.)  Then
$$
0 \qle T_\ell^{(m)} - T''_\ell \qle 
\phi_\ell(|r''_{k(m)-1} - r''_{k(m)+1}|) I_{\bar{R}_m}\,,
$$
so, using Lemma 5.3 in the second inequality below,
$$
\align
S_\ell \ :=\ \sum_{m=1}^\infty U_m 
&\qle \sum_m  \phi^2_\ell(|r''_{k(m)-1} - r''_{k(m)+1}|)I_{\bar{R}_m}\\
&\qequals \sum_{k=2}^{K''-1} \phi^2_\ell (|r''_{k-1} - r''_{k+1}|)\\
&\qle \sum_{k=2}^{K''-1} 2^{2\a} (\phi_\ell^2(|r''_{k-1}-r''_{k}|) + \phi_\ell^2(|r''_{k}-r''_{k+1}|) )\\
&\qle 2^{2\a+1}\sum_{k=1}^{K''-1} \phi^2_\ell (|r''_{k} - r''_{k+1}|)
\qequals 2^{2\a+1}S''_\ell.\tag\eqah\\
\endalign
$$
Combining (\eqaf), (\eqag), and (\eqah) gives 
$\var T''_\ell \le 2^{2\a+1}ES_\ell''$.

\enddemo

{\sl Step 2. For some constant $\cb$, $ES''_\ell \le \cb\ell$ (for 
$\ell > 1$).
In fact, with ${\ke} = 1/(4\a+2)$ and for some constants $\cs$ and $\cb$,}
$$
P[S''_\ell > x] \qle \cb \exp(-\cs x^{{\ke}})\text{\ for all } x\ge \cb\ell.\tag{\eqai}
$$

As this step is the heart of the proof, we begin by describing the overall
structure of the argument. A main ingredient is Lemma 3.3 below, which gives
a large deviation bound  for $T''_\ell$, obtained by constructing a suboptimal
path for the cost function $\phi_\ell$. Such arguments do not directly yield
bounds such as (3.10) for $S''_\ell$ because the definition of $S''_\ell$
involves replacing $\phi_\ell$ by $\phi_\ell^2$ while still using the 
links $L''_j$ that are optimal for $\phi_\ell$. So we separate the links
$L''_j$ into short and long ones and correspondingly write 
$S''_\ell = S_1 + \tilde{S}$ (with $\tilde{S}$ further decomposed as
$S_2 + S_3$). The tail of $S_1$ is directly estimated by that of $T''_\ell$,
but the analysis of $\tilde{S}$ requires more work. We will choose an
appropriately small $\e$, relate the path $r''$ to a kind
of path formed from $\e$-boxes and then control the tail of $\tilde{S}$
by a combination of percolation and lattice animal estimates for the
path formed from $\e$-boxes. Now, to work. 

We will call any finite sequence of distinct $\eta$-boxes an ``$\eta$-box path'' if the first
box contains the origin and the boxes are sequentially adjacent;  the path's ``length''
will refer to the number of boxes on the path.  
We call an $\eta$-box ``occupied'' if it contains a Poisson particle.
Pick $0<\e\le 1$ small enough so that, (i)
as in the proof of Lemma 3 of [HoN1], the events
$$
F_x \qequals \{\text{$\exists$ an $\e$-box path of length $m\ge x$ with at least $m/2d$
               occupied boxes}\}
$$
satisfy $P F_x \le (1-e^{-1})^{-1} e^{-x}$, and (ii) $17\e\sqrt{d}$ is
strictly less than the critical radius $R_c^*$ for continuum percolation
(discussed just before Conjecture 1 in Subsection 1.2).
The strict positivity of $R_c^*$ can be shown by standard arguments --- see, 
e.g., Theorem 3.2 of [MR].
(We remark that for any $\ell$, by the construction of 
$Q_\ell$ an $\e$-box is occupied (by a Poisson particle
in $Q$) if and only if it contains a particle in $Q_\ell$.)
Consider the $\e$-box path $\b=(\b_1,\dots,\b_{\tM(\x)})$ from $0$ to $\x$
constructed as
follows: $\b_1$ is the $\e$-box that contains $0$; if $\overline{r''}$
does not end inside of $\b_k$, $\b_{k+1}$
is the (\as\ adjacent) $\e$-box that $\overline{r''}$ enters when it last exits $\b_k$.
Here $\tM(\x)$ is the random number of boxes along this box path.
It follows as in the proof of Lemma 3 of [HoN1] that, for large $x$:
$$
T''_\ell\qge \frac{\phi_\ell(\e) x}{3d} \qequals \frac{\e^\a x}{3d}
\text{\ on } F_x^c\cap\{\tM(\x)\ge x\},
$$
(the equality above holds since $\e \le 1 \le h_\ell$) and hence
$$
P[\tM(\x)\ge x] \qle PF_x + P\Big[T''_\ell \ge \frac {\e^\a x}{3d}\Big]
\qle (1-e^{-1})\exp(-x) + P\Big[T''_\ell \ge \frac {\e^\a x}{3d}\Big].\tag{\eqaj}
$$

The $\e$-box path $\b$  covers
the midpoint of any sufficiently long link in $\overline{r''}$.
To see this,
let $\overline{ab} = \overline{r''_kr''_{k+1}}$ be any link in $\overline{r''}$ and
let $c$ be its midpoint.  Suppose $\b_{i^*}$ is the last
$\e$-box along $\b$ that touches either $\overline{ac}$ or any link that precedes
$\overline{ab}$ on $\overline{r''}$.  If $i^* = \tM(\x)$, put $\rho=\b_{i^*}$; otherwise
put $\rho=\b_{i^*}\cup \b_{i^*+1}$.  Then $\rho$ touches $c^*$ and $c^{**}$ satisfying
at least one of the following:
$$
\align
&c^*\in\overline{ac}\text{ and } c^{**}\in\overline{cb},\tag{\eqao}\\
&c^{*}\in L^*\text{ and } c^{**}\in \overline{cb}\text{ where $L^*$ is a link on $\overline{r''}$ preceding $\overline{ab}$,}\tag{\eqap}\\
&c^*\in \overline{ac}\text{ and } c^{**}\in L^{**}\text{ where $L^{**}$ is a link on $\overline{r''}$ succeeding $\overline{ab}$,}\tag{\eqaq}\\ 
&c^*\in L^*\text{ and }c^{**}\in L^{**}\text{ with $L^*$ and $L^{**}$ 
as in (\eqap) and (\eqaq).}\tag{\eqar}\\
\endalign
$$
Now (\eqao) implies that $c\in\rho$.  On the other hand, by the 
No Doubling Back Proposition of [Ho] (stated below as Lemma 5.5),  
(\eqap) implies
$$
\frac1 2 |a-b| \qequals |a-c| \qle |a-c^{**}| \qle 16|c^*-c^{**}| \qle 16\e\sqrt{d+3},
$$
while (\eqaq) similarly implies
$$
\frac1 2 |a-b| \qequals |c-b| \qle |c^{*}-b| \qle 16|c^*-c^{**}| \qle 16\e\sqrt{d+3},
$$
and (\eqar) implies
$$
\align
|a-b|&\qle|\text{ending point of }L^* - \text{starting point of }L^{**}|\\
& \qle 33|c^*-c^{**}|\qle 33\e\sqrt{d+3}.\\
\endalign
$$
It follows that $c\in\rho$ provided $|a-b|>33\e\sqrt{d+3}$. 

Choose $\l$ to be an odd integral multiple of $\e$ (so the $\e$-boxes
are nested in the $\l$-boxes) with $\l$ large enough that the probability that
any fixed $\l$-box contains no Poisson particle (equivalently, no $Q_\ell$ particle)
is below the critical probability for site
percolation on the nearest neighbor ${\Bbb Z}^d$ lattice.  
If the midpoint of a link $L''_k$ is touched by the $\e$-box path $\beta$, then a.s.{} it is
touched by only one of the $\e$-boxes on $\beta$; 
let $\nu(L''_k)$ denote the $\l$-box that contains
this $\e$-box.  
(If the midpoint of $L''_k$ is not so touched, $\nu(L''_k)$ is undefined.) For
any $\l$-box $\nu$, let $|{\Cal C}_\nu|$ denote the size 
(i.e., the cardinality) of the nearest-neighbor
cluster ${\Cal C}_\nu$ of unoccupied $\l$-boxes at $\nu$.  
The quantity $y_0$ in (\eqay) below will be specified later but
depends only on $d$.
We choose $h_0$ sufficiently large such 
that $|L''_k|>h_0$ implies: 
$$
\gather
\text{$\nu(L''_k)$ is defined;}\tag \eqaw \\
\text{if $L''_j\not=L''_k$ is another link with $|L''_j|>h_0$ then $\nu(L''_j)\not= \nu(L''_k)$; and}\tag\eqax\\
\text{$\nu(L''_k)$ is unoccupied, moreover $|{\Cal C}_{\nu(L''_k)}|\ge y_0^{1/(2\a)} |L''_k|$.}\tag\eqay\\
\endgather
$$
We can ensure (\eqaw) by the 
preceding discussion and (\eqax) also follows easily 
for large $h_0$ from the 
No Doubling Back Proposition (Lemma 5.5).  Since the interior of the region 
$\W_{\phi_\ell}(r''_k, r''_{k+1})$ 
contains no $Q_\ell$ particles,
Lemma 5.4 furnishes (\eqay) for $h_0$ sufficiently large (depending on $y_0$).

We split $S''_\ell$ into three pieces as follows:
$$
\align
S''_\ell &\qequals S_1+S_2+S_3,\text{\ where}\\
S_1 &\qequals \sum_{k\ :\ |L''_k|\le h_0} \phi_\ell^2(|L''_k|),\\
S_2 &\qequals I_{\{\tM(\x)\ge x\}} \sum_{k\ :\ |L''_k|>h_0} \phi_\ell^2(|L''_k|), \text{\ and}\\
S_3 &\qequals I_{\{\tM(\x)< x\}} \sum_{k\ :\ |L''_k|>h_0} \phi_\ell^2(|L''_k|).\\
\endalign
$$
Now $S_1 \le h_0^\a T''_\ell$, so
$$
P[S_1 > x] \qle \cb \exp(-\cs x^{{\kf}})\text{\ for all } x\ge \cb\ell\tag{\eqak}
$$
will follow with ${\kf} = \min(1,d/\a)$ from:
\proclaim{Lemma 3.3} There exist constants $\cs$ and $\cb$
such that, for $T = T_\ell$,  $T=T'_\ell$, or $T=T''_\ell$, 
$P[T>x] \le \cb \exp(-\cs x^{{\kf}})$ for $x \ge \cb \ell$.
\endproclaim
\demo{Proof}  We first prove (in detail) the 
case $T=T'_\ell$. For $a\in \rd$ and $t\ge0$,
let
$$\align
{\Cal D}_t(a) &\qequals \{a+b\in\rd:0\le b_1\le t; 0\le \sigma_i b_i\le b_1\text{ for }2\le i\le d\},
\text{\ where}\cr
\sigma_i &\qequals -1\text{ if }a_i\ge0,\ 1\text{ otherwise.}\cr
\endalign
$$
Then the $d$-dimensional volume of 
${\Cal D}_t(a)$ is $\int_0^t s^{d-1}\,ds = t^d/d$.
Also, for $b\in {\Cal D}_t(a)$ we have
$$
\max_{2\le i\le d} |b_i| \qle \max(t, \max_{2\le i\le d} |a_i|). \tag{\eqal}
$$
Let $q_0 = 0$ and define $q_n$ 
and $R_n$ inductively for $n\ge1$ (See Figure 1 for the picture when $d=2$)
by the relation
$$
\tR_n \qequals \inf\{t>0:\text{ there exists a Poisson particle }q_n\not=q_{n-1}\text{ in }{\Cal D}_t(q_{n-1})\},
$$
and let $\tR^*_n = \max_{1\le m\le n} \tR_m$.

\vskip .5truein
\beginpicture
\setcoordinatesystem units <.4in,.4in>
\setplotarea x from -1 to 11, y from -5 to 5
\plotheading {Figure 1.}

\axis bottom /
\axis top /
\axis left /
\axis right /

\put {$\cdot$} at     0.59     4.45
\put {$\cdot$} at     3.74    -1.84
\put {$\cdot$} at     2.87     4.71
\put {$\cdot$} at     5.22    -3.75
\put {$\cdot$} at     2.62    -3.04
\put {$\cdot$} at     1.98     4.70
\put {$\cdot$} at     0.37    -4.99
\put {$\cdot$} at     4.39    -1.64
\put {$\cdot$} at     3.17     1.68
\put {$\cdot$} at     6.74    -0.26
\put {$\cdot$} at     4.95    -3.53
\put {$\cdot$} at     7.38    -2.07
\put {$\cdot$} at     9.27     3.22
\put {$\cdot$} at     9.15     3.00
\put {$\cdot$} at     2.19     4.88
\put {$\cdot$} at     1.61    -4.58
\put {$\cdot$} at     0.41    -2.42
\put {$\cdot$} at     8.89    -4.64
\put {$\cdot$} at     7.38     1.19
\put {$\cdot$} at     9.62     1.37
\put {$\cdot$} at     2.34    -1.33
\put {$\cdot$} at     7.80     3.26
\put {$\cdot$} at     5.28     1.12
\put {$\cdot$} at     1.10    -1.16
\put {$\cdot$} at     5.37    -0.80
\put {$\cdot$} at     1.85     2.16
\put {$\cdot$} at     3.34    -0.80
\put {$\cdot$} at     1.81     0.16
\put {$\cdot$} at     9.55     1.84
\put {$\cdot$} at     0.22     2.26
\put {$\cdot$} at     5.78    -1.10
\put {$\cdot$} at     1.45    -4.16
\put {$\cdot$} at    -0.87    -3.88
\put {$\cdot$} at    10.69     3.53
\put {$\cdot$} at     3.91    -3.04
\put {$\cdot$} at     4.14     0.85
\put {$\cdot$} at     4.10    -4.30
\put {$\cdot$} at    -0.26    -2.89
\put {$\cdot$} at     7.53     4.88
\put {$\cdot$} at     5.08     2.18
\put {$\cdot$} at     8.21    -3.40
\put {$\cdot$} at     1.71    -4.81
\put {$\cdot$} at    10.69     0.48
\put {$\cdot$} at     1.25    -3.74
\put {$\cdot$} at     8.18    -3.77
\put {$\cdot$} at     7.96     3.87
\put {$\cdot$} at     5.10    -2.69
\put {$\cdot$} at     8.65    -2.68
\put {$\cdot$} at     2.10     0.89
\put {$\cdot$} at     1.52    -4.83
\put {$\cdot$} at     6.33     3.88
\put {$\cdot$} at     7.03     3.12
\put {$\cdot$} at     9.56    -0.81
\put {$\cdot$} at     9.91    -4.97
\put {$\cdot$} at     3.97     3.98
\put {$\cdot$} at     7.06    -0.40
\put {$\cdot$} at     4.49     0.22
\put {$\cdot$} at     4.94     3.26
\put {$\cdot$} at     1.48    -1.80
\put {$\cdot$} at    -0.77     3.07
\put {$\cdot$} at     3.34     4.16
\put {$\cdot$} at    -0.84    -0.87
\put {$\cdot$} at     3.58    -4.41
\put {$\cdot$} at     1.44    -2.26
\put {$\cdot$} at     7.92    -2.43
\put {$\cdot$} at     1.57     2.14
\put {$\cdot$} at     2.86    -0.12
\put {$\cdot$} at     2.49     3.32
\put {$\cdot$} at     9.61    -0.54
\put {$\cdot$} at     3.41    -2.11
\put {$\cdot$} at    -0.33    -1.98
\put {$\cdot$} at     3.92    -4.72
\put {$\cdot$} at     5.94     1.48
\put {$\cdot$} at     2.26    -3.55
\put {$\cdot$} at     3.84     1.39
\put {$\cdot$} at     5.35     1.92
\put {$\cdot$} at     5.20    -2.84
\put {$\cdot$} at     0.31    -3.05
\put {$\cdot$} at     8.55     4.22
\put {$\cdot$} at     6.69    -1.96
\put {$\cdot$} at    10.03     1.46
\put {$\cdot$} at     2.05    -3.21
\put {$\cdot$} at     8.17    -0.85
\put {$\cdot$} at     2.53     4.69
\put {$\cdot$} at    10.41    -1.39
\put {$\cdot$} at     1.35     4.03
\put {$\cdot$} at     5.64     4.73
\put {$\cdot$} at     8.15     3.11
\put {$\cdot$} at     6.30    -3.59
\put {$\cdot$} at     1.25     0.76
\put {$\cdot$} at     6.64    -4.72
\put {$\cdot$} at     1.25     0.08
\put {$\cdot$} at    10.41    -0.96
\put {$\cdot$} at     6.63     4.36
\put {$\cdot$} at     5.23     0.79
\put {$\cdot$} at     0.11     2.52
\put {$\cdot$} at     6.45     4.50
\put {$\cdot$} at    10.94    -4.11
\put {$\cdot$} at     0.65     2.48
\put {$\cdot$} at     0.68     2.78
\put {$\cdot$} at     8.89     2.78
\put {$\cdot$} at    -0.30     3.11
\put {$\cdot$} at     9.98    -2.96
\put {$\cdot$} at     4.14    -1.04
\put {$\cdot$} at     8.21    -1.73
\put {$\cdot$} at     9.96    -1.84
\put {$\cdot$} at     4.19     0.63
\put {$\cdot$} at     4.62     2.22
\put {$\cdot$} at     7.07    -3.52
\put {$\cdot$} at    10.94     4.72
\put {$\cdot$} at     3.56     2.42
\put {$\cdot$} at     6.74     1.81
\put {$\cdot$} at     8.64    -3.47
\put {$\cdot$} at    -0.68    -0.73
\put {$\cdot$} at     0.87     4.00
\put {$\cdot$} at     0.92    -4.05
\put {$\cdot$} at     2.25     3.33
\put {$\cdot$} at     4.33    -1.71
\put {$\cdot$} at     0.46     0.27
\put {$\cdot$} at     4.86    -4.51
\put {$\cdot$} at     3.44     1.69
\put {$\cdot$} at     2.66     0.13
\put {$\cdot$} at    -0.79    -3.38
\put {$\cdot$} at     1.82     3.90
\put {$\cdot$} at     5.83    -3.44
\put {$\cdot$} at     0.12     3.06
\put {$\cdot$} at    -0.64     2.33
\put {$\cdot$} at     0.80     3.46
\put {$\cdot$} at     9.88    -0.44
\put {$\cdot$} at    10.34    -0.38
\put {$\bullet$}  at     0.00     0.00
\put {0} at 0 .25
\multiput {$\cdot$} at 0 0 *100 .1 0 /
\put {$\bullet$}  at    10.00     0.00
\put {$\ell\hat e_1$} at 10.1 .25
\put {$\ast$} at     2.34    -1.33
\put {$q_1$} [r] at 2.24 -1.33
\plot     0.00     0.00      2.34     0.00 /
\plot     2.34     0.00      2.34    -2.34 /
\plot     2.34    -2.34      0.00     0.00 /
\put {$\ast$} at     3.34    -0.80
\put {$q_2$} [lt] at      3.44   -.90
\plot     2.34    -1.33      3.34    -1.33 /
\plot     3.34    -1.33      3.34    -0.33 /
\plot     3.34    -0.33      2.34    -1.33 /
\put {$\ast$} at     4.49     0.22
\put {$q_3$} [lb] at 4.59     0.32
\plot     3.34    -0.80      4.49    -0.80 /
\plot     4.49    -0.80      4.49     0.34 /
\plot     4.49     0.34      3.34    -0.80 /
\put {$\ast$} at     6.69    -1.96
\put {$q_4$} [t] at  6.69    -2.06
\plot     4.49     0.22      6.69     0.22 /
\plot     6.69     0.22      6.69    -1.99 /
\plot     6.69    -1.99      4.49     0.22 /
\put {$\ast$} at     8.17    -0.85
\put {$q_5$} [lt] at 8.27   -0.95
\plot     6.69    -1.96      8.17    -1.96 /
\plot     8.17    -1.96      8.17    -0.48 /
\plot     8.17    -0.48      6.69    -1.96 /
\put {$\ast$} at     9.56    -0.81
\put {$q_6$} [t] at  9.56    -0.91
\plot     8.17    -0.85      9.56    -0.85 /
\plot     9.56    -0.85      9.56     0.53 /
\plot     9.56     0.53      8.17    -0.85 /
\put {$\ast$} at    10.34    -0.38
\put {$q_7$} [l] at 10.44 -0.38
\plot     9.56    -0.81     10.34    -0.81 /
\plot    10.34    -0.81     10.34    -0.03 /
\plot    10.34    -0.03      9.56    -0.81 /

\endpicture
\vskip .5truein

Now $|q_{n-1} - q_n| \le \tR_n\sqrt{d}$ and it follows from (\eqal) that 
$|q_N - \x| \le \tR^*_N\sqrt{d}$ where
$$
N\qequals \min\{n:\tR_1+\dots+\tR_n\ge\ell\}.
$$
Hence
$$
\align
T'_\ell &\qle |q_0 - q_1|^\a +\dots+|q_{N-1}-q_N|^\a + |q_N - \x|^\a\cr
  &\qle (\tR_1\sqrt{d})^\a +\dots+(\tR_N\sqrt{d})^\a+(\tR^*_N\sqrt{d})^\a\cr
  &\qle 2d^{\a/2} (\tR_1^\a+\dots+\tR_N^\a).\cr
\endalign
$$
It follows that for any $n>0$
$$
\align
P[T'_\ell>x] &\qle P[2d^{\a/2} (\tR_1^\a+\dots+\tR_n^\a)>x]\ +\ P[n<N]\\
       &\qle P[2d^{\a/2} (\tR_1^\a+\dots+\tR_n^\a)>x]\ +\ P[\tR_1+\dots+\tR_n <\ell]\tag{\eqam}\\
\endalign
$$ 
Now the $\tR_i$, and 
hence the $\tR_i^\a$, are i.i.d.{}, with $P[\tR_i^\a>r] = P[\tR_i>r^{1/\a}]
= \exp(-\frac1 d r^{d/\a})$. 
Taking $n = \lceil cx\rceil$ in (\eqam), it follows from [N] that, 
for sufficiently small $c$, there exist $\cs$ and $\cb$ such that
$$
P[2d^{\a/2} (\tR_1^\a+\dots+\tR_{\lceil cx\rceil}^\a)>x] 
\qle  \cb  \exp(-\cs x^{{\kf}})\text{\ for all }x.
$$
Also, for this choice of $c$, it follows 
from elementary large deviation results for i.i.d. random
variables (see, e.g., Sec. 1.9 of [D]) that,
for possibly larger $\cb$ and smaller $\cs$, we have:
$$
P[\tR_1+\dots+\tR_{\lceil cx\rceil} <\ell] 
\qle  \cb \exp(-\cs x)\text{\ for } x\ge \cb\ell.
$$
The lemma therefore follows for $T'_\ell$.

This extends easily to $T=T_\ell$ (with the same exponent ${\kf}$) 
by applying the first part of Lemma 3.1. 
To apply the $T'_\ell$ result to $T''_\ell$, note that the fact that
$\phi_\ell(t)\le \phi_\infty(t)=t^\a$ is helpful, so the only difficulty is
that the sequence of Poisson particles $q_1,\dots,q_N$ constructed above
are not necessarily in $Q_\ell$. 
However, there is always a $Q_\ell$ particle $\tilde q_i$
within
a distance $(\e/3^\ell)\sqrt{d} \le \sqrt{d}$ of each $q_i$ constructed above.  It is not
hard to see that the sequence $(\tilde q_i)$ 
produces a path whose passage time has a distribution with the requisite tail,
again with the same exponent ${\kf}$.\qed
\enddemo

We bound the tail of $S_2$ by the simple estimate
$$
\align
P[S_2>x]&\qle P[\tM(\x)>x]\\
        &\qle (1-e^{-1})\exp(-x) + \cb\exp\bigg(-\cs\Big(\frac {\e^\a x}{3d}\Big)^{{\kf}}\bigg)
\text{\ for $\frac {\e^\a x}{3d}\ge \cb\ell$}\\
&\qle \cb \exp(-\cs x^{{\kf}})\text{\ for all } x\ge \cb\ell.\tag{\eqan} 
\endalign
$$
Here we use (\eqaj) and Lemma 3.3; the final inequality holds 
for possibly larger $\cb$ and smaller $\cs$ since $\kf \le 1$.

Finally, we bound the tail of $S_3$.
Let $\xi''$ denote the
collection of $\l$-boxes that contain an $\e$-box on $\b$.   If 
$\Xi_0 = \{\text{all $\zd$ lattice animals containing the origin}\}$, then
$\xi''\in\Xi_0$ in the sense that the sites in $\zd$ associated
with the boxes in $\xi''$ form a lattice animal containing the origin.
Then, using (\eqaw), (\eqax), and (\eqay),
$$
\align
S_3 &\qle I_{\{\tM(\x)< x\}}\sum_{k\ :\ |L''_k|>h_0} |L''_k|^{2\a}\\
    &\qle I_{\{\tM(\x)< x\}} \sum_{k\ :\ |L''_k|>h_0} y_0^{-1}|{\Cal C}_{\nu(L''_k)}|^{2\a}\\
    &\qle I_{\{\tM(\x)< x\}} \sum_{\nu\in\xi''} y_0^{-1}|{\Cal C}_\nu|^{2\a},\\
\endalign
$$
and hence, using that  $|\xi''|$, the number of sites (boxes) in $\xi''$, cannot
exceed $\tM(\x)$, we have for any $\g < 1$,
$$
\align
\{S_3 > x\} \ \subset\ &\Big\{|\xi''| < x^\g\ \text{ and } \sum_{\nu\in\xi''} |{\Cal C}_{\nu}|^{2\a} > y_0 x\Big\}\\
  &\cup\ \Big\{x^\g \le |\xi''| \le x \text{ and } \sum_{\nu\in\xi''} |{\Cal C}_{\nu}|^{2\a} > y_0 x\Big\}\\
  \ \subset\ &\Big\{\exists \nu\in [-x^\g, x^\g]^d\cap\zd \text{ with } |{\Cal C}_\nu| > y_0^{1/(2\a)} x^{(1-\g)/(2\a)}\Big\}\\
             &\cup\ \Big\{\text{$\exists\xi\in\Xi_0$  with $|\xi| \ge x^\g$  and  $\frac{1}{|\xi|}\sum_{\nu\in\xi}|{\Cal C}_\nu|^{2\a} > y_0$}\Big\}.\tag\eqas\\
\endalign
$$
But, for some constant $b>0$, $P[|{\Cal C}_\nu| > x] \le \exp(-bx)$ for all $x$
(see, e.g., [Gr]), so
$$
\align
P\Big[\exists \nu\in [-x^\g, x^\g]^d\cap\zd \text{ with } &|{\Cal C}_\nu| > y_0^{1/(2\a)} x^{(1-\g)/(2\a)}\Big]\\ 
&\qle (2x^\g+1)^d \exp(-by_0^{1/(2\a)}x^{(1-\g)/(2\a)}).\tag\eqat\\
\endalign
$$
By Theorem 5 of [HoN2], provided $y_0$ 
is sufficiently large (depending only on
$d$ and the distribution of 
the $|{\Cal C}_\nu|$, which in turn depends only on $d$), 
we also have
for some $a>0$ and a possibly smaller $b$:
$$
P\Big[\text{$\exists\xi\in\Xi_0$ with $|\xi| \ge x^\g$ and $\frac{1}{|\xi|}\sum_{\nu\in\xi}|{\Cal C}_\nu|^{2\a} > y_0$}\Big]
 \qle a \exp (-bx^{\g/(2\a+2)}).\tag\eqau
$$
The exponents $(1-\g)/(2\a)$ in (\eqat) and $\g/(2\a+2)$ in (\eqau) are
both made equal to ${\ke}$ by taking $\g = (\a+1)/(2\a+1)$.  For this choice of $\g$, 
combining (\eqas), (\eqat) and (\eqau) gives that
$$
P[S_3 > x] \qle \cb \exp(-\cs x^{{\ke}})\text{\ for all } x\tag\eqav
$$
for possibly some larger $\cb$ and smaller $\cs$. Noting
that ${\ke} < {\kf}$, combining (\eqak), (\eqan), and (\eqav) yields that
$$
P[S''_\ell>3x]\qle \cb \exp(-\cs x^{{\ke}})\text{\ for all } x\ge \cb\ell,
$$
which proves (\eqai) for possibly larger $\cb$ and smaller $\cs$.
Step 2 is completed as follows:
$$
\align
ES''_\ell &\qequals \int_0^\infty P[S''_\ell > x]\,dx\\
          &\qle \cb\ell + \int_{\cb\ell}^\infty \cb\exp(-\cs x^{{\ke}})\,dx\\
          &\qequals \cb\ell + o(\ell)\text{\ as $\ell\to\infty$}.\\
\endalign
$$

{\sl Step 3}.  $\var T_\ell \le \cb \ell$ {\sl\ for\ } $\ell>1$.
Steps 1 and 2 show that, for appropriate $\cb$, $\var T''_\ell < \cb \ell$ for 
$\ell>1$.  Now
$$
\align
\std T_\ell &\qle \std T''_\ell + \std |T'_\ell - T''_\ell| + \std |T'_\ell - T_\ell|\\
&\qle \cb^{1/2}\ell^{1/2} + \std |T'_\ell - T''_\ell| + \std |T'_\ell - T_\ell|.\\
\endalign
$$
It follows from (\eqac) 
that $\std |T'_\ell - T_\ell|$ is bounded in $\ell$.  On the
other hand, since $0 \le T'_\ell, T''_\ell \le \ell^\a$, we have
$|T'_\ell - T''_\ell| \le \ell^\a I_{\{T'_\ell \not= T''_\ell\}}$ 
so, assuming (\eqad) holds,
$$
\var |T'_\ell - T''_\ell| \qle E(|T'_\ell - T''_\ell|^2)
\qle \ell^{2\a} P[T'_\ell \not= T''_\ell] 
\qequals o(\ell) \text{\ as $\ell\to\infty$},
$$
yielding that, for possibly larger $\cb$, $\var T_\ell < \cb \ell$ for 
$\ell>1$.  In view of Lemma 3.2, (2.4) will be proved once we complete the:

\demo{Proof of (\eqad)} For an $a>1$ (to be chosen momentarily), let 
$B(a\ell) = [-a\ell,a\ell]^d$. If $B$ is any cube containing $0$ and $\x$
such that $\overline{r'}\subset B$, $\overline{r''}\subset B$,
$Q\cap B = Q_\ell\cap B$, and 
no link on $\overline{r''}$ exceeds $h_\ell$ in length,
then $T'_\ell = T''_\ell$.  Hence
$$
\aligned
P[T'_\ell \not= T''_\ell] &\qle P[\overline{r'} \not\subset B(a\ell)] + 
                                P[\overline{r''} \not\subset B(a\ell)]\\
                          &\ + P[\exists \text{\ an 
 $(\e/3^{\lfloor\ell\rfloor})$-box 
 touching $B(a\ell)$ with two or more Poisson particles}]\\
          &\ + P[\exists \text{\ a $\l$-box $\nu$ touching $B(a\ell)$ with
                            $|{\Cal C}_\nu| \ge y_0^{1/(2\a)}h_\ell$}]\, ,\\
\endaligned\tag\eqbw
$$
where we used (\eqaw) and (\eqay).
First, we bound the term $P[\overline{r''}\not\subset B(a\ell)]$.
If $\overline{r''}\not\subset B(a\ell)$, then either
$\beta\not\subset B(a\ell/2)$ or else $\beta\subset B(a\ell/2)$
and for some $(r''_{i_1},r''_{i_1+1},\dots,
r''_{i_2},\dots,r''_{i_2})$
we have that $\overline{r''_{i_1}\,r''_{i_1+1}}$ exits an
$\e$-box $\beta_k$ on $\beta$, $r''_{i_2}\not\in B(a\ell)$, 
and $\overline{r''_{i_3-1}\,r''_{i_3}}$ re-enters $\beta_k$.
By the No Doubling Back Proposition (Lemma 5.5), in the latter
case we must have that $r''_{i_1+1}$ and $r''_{i_3-1}$ are within
Euclidean distance $16\e\sqrt{d}$ of $\beta_k\subset B(a\ell/2)$
and also that
$|r''_{i_1+1}-r''_{i_3-1}| \le 33\e\sqrt{d}$. It follows 
also (since $r''$ is minimizing) that $|r''_i-r''_{i+1}|\le33\e\sqrt{d}$
for $i_1< i<i_3-1$. These together would imply that there is a cluster of
overlapping balls of radius $17\e\sqrt{d}$ centered at
particle locations in $Q$ touching both $B(a\ell/2)$ and 
$B(a\ell)^c$.  Since $17\e\sqrt{d}$ is less than the critical
continuum percolation radius $R_c^*$, 
this latter event occurs with
probability bounded by $\cb\exp(-\cs\ell)$ --- a consequence
of Theorem 3.5 and Lemma 3.3 of [MR] (here $\cs$ and $\cb$
depend on $d$, $\e$ and $a$). 
It follows that
$$
P[\overline{r''}\not\subset B(a\ell)] \qle
P[\beta\not\subset B(a\ell/2)]\ +\ \cb\exp(-\cs\ell).
$$
We take $\cb$ as in the rightmost expression 
of (3.22) and then for sufficiently large $a$,
we have from the definitions of $\beta$ and ${\tM(\x)}$
that $\{\beta\not\subset B(a\ell/2)\}\subset\{\tM(\x)>\cb\ell\}$
so, as in (3.22), $P[\beta\not\subset B(a\ell/2)]\le\cb\exp(-\cs\ell^\kf)$.
Since $\kf\le1$, this yields 
$P[\overline{r''}\not\subset B(a\ell)]\le\cb\exp(-\cs\ell^\kf)$.
The first term on the right side 
of (\eqbw) may be similarly bounded for a possibly
larger $a$.  

With $a$ now
fixed, there are $O(\ell^d3^{\ell d})$ $(\e/3^{\lfloor\ell\rfloor})$-boxes 
and $O(\ell^d)$ $\l$-boxes touching
$B(a\ell)$.  Since the probability 
that any particular $(\e/3^{\lfloor\ell\rfloor})$-box
has two or more Poisson particles in it is 
bounded by $(\e/3^{\lfloor\ell\rfloor})^{2d}$,
the third term on the right side of (\eqbw) is of 
order $\ell^d 3^{-\ell d}\le \cb\exp(-\ell)$ for possibly larger
$\cb$.  Finally, by our earlier choice of $\l$,
the probability that any particular $\l$-box $\nu$ has
$|{\Cal C}_\nu| \ge y_0^{1/(2\a)}h_\ell$ is bounded by $\exp(-by_0^{1/(2\a)}h_\ell)$ yielding
that the fourth term in (\eqbw) is bounded by $\cb\exp(-\cs\ell^{1/(2\a)})$ for possibly 
larger $\cb$ and smaller $\cs$ since $h_1 > 0$.  
Collectively this proves (\eqad) since ${\kc} = 1/(2\a) < {\kf} \le 1$.\qed
\enddemo

This completes the proof of (2.4). We finish the proof of Theorem 2.1 with:

{\sl Step 4.  Proof of (2.5).}
Our strategy here is to invoke Lemma 5.6 for large $\ell$, using $\field_m$,
$\Delta_m$, and $U_m$ from the previous section, i.e.,
$$
\Delta_m \qequals E[T''_\ell | \field_m] - E[T''_\ell| \field_{m-1}]
\text{\ and }
U_m \qequals (T^{(m)}_\ell - T''_\ell)^2.
$$
We also therefore take $S = S_\ell$ as given in (\eqah).  We presently
show that the hypotheses of the lemma are satisfied for appropriate $x_0$, $c$, and
$\gamma$.

First, we observe that 
$0 \le T_\ell^{(m)} - T''_\ell \le 2^\a h_\ell^\a$.
The first inequality is trivial and the second follows from Lemma 5.3. Since
$T^{(m)}_\ell$ is independent of $\field(B_m)$ we see that
$E[T^{(m)}_\ell | \field_m] = E[T^{(m)}_\ell | \field_{m-1}]$.  It follows that
$|\Delta_m| \le 2^\a h_\ell^\a$. We therefore take $c = 2^\a h_\ell^\a$ 
in Lemma 5.6.

Next, we verify that 
$E[ \Delta_m^2 | \field_{m-1}] \le E[U_m | \field_{m-1}]$ as
follows:
$$
\aligned
E[ \Delta_m^2 |\field_{m-1}] 
&\qequals E[(E[T''_\ell | \field_m] - E[T''_\ell| \field_{m-1}])^2 |\field_{m-1}]\\
&\qle E[(E[T''_\ell | \field_m] - E[T^{(m)}_\ell| \field_m])^2 |\field_{m-1}]\\
&\qequals E[ (E[T''_\ell - T^{(m)}_\ell | \field_m])^2 | \field_{m-1}]\\
&\qle     E[ E[(T''_\ell - T^{(m)}_\ell)^2 | \field_m] | \field_{m-1}]\\
&\qequals E[(T''_\ell - T^{(m)}_\ell)^2 | \field_{m-1}] \qequals E[U_m | \field_{m-1}].\\
\endaligned
$$
The first inequality uses (\eqbx) with 
${\Cal G} = \field_{m-1}$, $X = E[T''_\ell | \field_m]$,
and
$Y = E[T^{(m)}_\ell | \field_{m-1}] = E[T^{(m)}_\ell | \field_m]$.  The second 
inequality follows from the conditional Jensen's inequality.

Additionally, by (\eqah) and (\eqai) and with ${\kb} = 1/(4\a+3) < {\ke}
= 1/(4\a+2)$, we get
$$
\aligned
P[S>x] &\qle P[S''_\ell > 2^{-(2\a+1)}x]\\
&\qle \cb \exp\big(-\cs (2^{-(2\a+1)}x)^{{\ke}}\big) \text{\ for $x\ge 2^{2\a+1}\cb\ell$}\\
&\qle \cb\exp(-x^{\kb}) \text{\ for $x\ge 2^{2\a+1}\cb\ell$}, \\
\endaligned\tag\eqbnew
$$
where the last inequality holds for a possibly larger
$\cb$.  This gives (\eqbj) with 
$\g = {\kb}$ and $x_0 = 2^{2\a+1}\cb\ell$.  

Finally, we must have $x_0 \ge c^2 \ge 1$.  The first inequality holds if
$h_\ell \le (2\cb\ell)^{1/(2\a)}$.  Recalling that $h_\ell = \max(h_0, h_1\ell^{1/(2\a)})$
where $h_0$ has already been specified, we take $h_1 = (2\cb)^{1/(2\a)}$.  
We will
then have $x_0\ge c^2$ for $\ell$ large enough that $h_\ell = h_1\ell^{1/(2\a)}$.  
The second inequality ($c\ge1$) is equivalent
to $h_\ell \ge 1/2$ which holds since $h_0\ge1$.

Lemma 5.6 implies that there are constants $\cs$  and $\cb$ such that,
for $\ell$ large enough that $h_\ell = h_1\ell^{1/(2\a)}$, 
$$
\prob[\,|T''_\ell-ET''_\ell| > x\sqrt{\ell}] \qle \cb\exp(-\cs x)
\text{\ for }x\le \cs\ell^{{\kb}},
$$
which can be made to hold for all $\ell$ by increasing $\cb$.
Now 
$$
|T''_\ell - T_\ell| \qle |T''_\ell - T'_\ell| + |T'_\ell - T_\ell|
\qle \ell^\a I_{\{T''_\ell \not= T'_\ell\}} + |T'_\ell - T_\ell|,
$$
so it follows from Lemma 3.1 that $|ET''_\ell - ET_\ell|$ is bounded by some 
constant $\tilde b$.
Also, using that
$$
|T_\ell - ET_\ell| \qle |T_\ell - T'_\ell| + |T'_\ell - T''_\ell| + |T''_\ell - ET''_\ell|
+ |ET''_\ell - ET_\ell|,
$$
we get, for $\ell>1$ and ${\tilde b}\le x\le \cs \ell^\kb$, that
$$
\align
P[|T_\ell - ET_\ell| > 3x\sqrt{\ell}] \qle &P[|T_\ell - T'_\ell|>x\sqrt{\ell}]
                                             + P[T'_\ell \not= T''_\ell]\\
                                           &\qquad  + P[|T''_\ell - ET''_\ell| >x\sqrt{\ell}]\\
                                      \qle &\cb \exp(-\cs (x\sqrt{\ell})^\kd)
                                            + \cb \exp(-\cs \ell^\kc)\\
                                            &\qquad + \cb \exp(-\cs x).\tag\eqbr\\
\endalign
$$
On the one hand, (\eqbr) produces for appropriate $\cs$ and $\cb$ and
for $\ell>1$ and ${\tilde b}\le x\le \cs \ell^\kb$,
$$
\align
P[|T_\ell - ET_\ell| > 3x\sqrt{\ell}] \qle &\cb\exp(-\cs x^\kd)
                                            + \cb\exp(-\cs x^{\kc/\kb})
                                            + \cb \exp(-\cs x)\\
                                      \qle &\cb\exp(-\cs x^\kf),\\
\endalign
$$
with the last inequality holding for possibly larger $\cb$ since 
$\kf = \min(1, \kd)$ and $\kc/\kb > 1$.  By possibly increasing $\cb$ still further
and decreasing $\cs$ we can ensure that
$$
P[|T_\ell - ET_\ell| > x\sqrt{\ell}] \qle \cb\exp(-\cs x^\kf)
\text{ for all $\ell$ and $x\le \cs\ell^\kb$,}
$$
proving (2.5).\qed

\head{4. Proof of Theorems 2.2 and 2.3}\endhead 

Our plan is to show that
$ET_\ell$ exhibits the following sort of weak superadditivity:
\proclaim{Lemma 4.1} For some constant $\cb\in (0,\infty)$ 
we have
$$
ET_{2\ell} \qge 2ET_\ell - \cb\sqrt{\ell} (\log \ell)^{1/\kappa_1}
\text{\quad for all large $\ell$}. \tag{\eqca}
$$
\endproclaim

Before proving this lemma, which constitutes the bulk of this
section, we show how this gives Theorems 2.2 and 2.3.  First, we need the 
following easy lemma; it will be applied with $a(\ell)=ET_\ell$
and $g(l)=\cb\sqrt{\ell} (\log \ell)^{1/\kappa_1}$.

\proclaim{Lemma 4.2} Suppose the functions $a:\rplus\to \reals$ 
and $g:\rplus\to \rplus$ satisfy
the following conditions: $a(\ell)/\ell\to\nu\in \reals$, 
$g(\ell)/\ell\to0$ as $\ell\to\infty$, $a(2\ell) \ge 2a(\ell) - g(\ell)$,
and $\psi\equiv\limsup_{\ell\to\infty} g(2\ell)/g(\ell) < 2$.
Then, for any $c>1/(2-\psi)$, 
$a(\ell) \le \nu\ell + cg(\ell)$ for all large $\ell$.
\endproclaim
\demo{Proof}  It is easily verified that, for $c>1/(2-\psi)$,
$\tilde a(\ell)\equiv a(\ell)-cg(\ell)$ satisfies 
$\tilde a(2\ell) \ge 2\tilde a(\ell)$ for all large $\ell$.  
Iterating this $n$ times
yields $\tilde a(2^n\ell) \ge 2^n\tilde a(\ell)$ or 
$\tilde a(2^n\ell)/(2^n\ell)\ge \tilde a(\ell)/\ell$.  
Under our hypotheses on $a$ and $g$,
$\tilde a(x)/x\to\nu$ as $x\to\infty$, so
letting $n\to\infty$ shows that $\tilde a(\ell) /\ell \le \nu$
for all large $\ell$.\qed
\enddemo

\demo{Proof of Theorems 2.2 and 2.3}  Based 
on general subadditivity considerations, we have (see [HoN1]) that 
$$
0 \ <\ \mu\ \equiv\ \inf_{\ell>0} \frac{ET_\ell}{\ell}\ <\ \infty
\text{\quad and\quad}
\lim_{\ell\to\infty} \frac{T_\ell}{\ell} \qequals \mu
\text{\ (\as\ and in $L^1$)}. \tag{\eqcb}
$$
Taking
$a(\ell) = ET_\ell$ and
$g(\ell) = \cb\sqrt{\ell}(\log\ell)^{1/\kappa_1}$ 
in Lemma 4.2 (so that $\limsup_\ell g(2\ell)/g(\ell)$ $= \sqrt{2}<2$),
we get that, for appropriate $\cb$,
$$
\mu\ell \qle ET_\ell \qle \mu\ell + \cb \sqrt{\ell}(\log\ell)^{1/\kappa_1}
\text{\ for large $\ell$}. \tag{\eqcc}
$$

The second part of Theorem 2.1 then immediately implies that
$$
P[|T_\ell - \mu\ell| > 2x\sqrt\ell]\qle \cb\exp(-\cs x^\kf)
\quad\text{for $\cb(\log\ell)^{1/\kf} \le x \le \cs\ell^\kb$}.
$$
Substituting $\l = 2x\sqrt\ell$ yields (2.7) for large $\ell$, with
this latter restriction lifted by adjusting $\cs$ and $\cb$, which
proves Theorem 2.2.
On the other hand, substituting $x = \frac1 2 (\log\ell)^{(1+\e)/\ka}$, where $\e>0$, yields
$$
P[|T_\ell - \mu\ell| > \sqrt\ell(\log\ell)^{(1+\e)/\ka}]\qle 
\cb\ell^{-\cs{(\log\ell)^\e}}\quad\text{for large $\ell$}.
$$
This and the Borel-Cantelli Lemma together imply that, \as, the
event $\{\big|T(0,w) - \mu|w|\big| > \sqrt{|w|}(\log|w|)^{(1+\e)/\ka}\}$
occurs for only finitely many $w\in\zd$.  Theorem 2.3 follows from this 
together with
an application of Lemma 5.2 and the Borel-Cantelli Lemma. Further details
are left to the reader.\qed
\enddemo

\demo{Proof of Lemma 4.1} Fix $\gamma$ 
with $0<\a\g < 1/2$. Define the event
$$
\F_\ell \ \equiv\ \{\text{there 
exists an $x\in\rd$ 
with $|x-\ell\hat e_1| \le 3\ell$ and $|q(x)-x|\ge \ell^\gamma$}\}.
$$

Next, take $x_1 = \ell \hat e_1$ and pick $x_2,\dots,x_{n(\ell)}$
on $\partial \B(0,\ell)$, the Euclidean sphere of radius $\ell$
centered at the origin, so that every $x\in \partial \B(0,\ell)$
is within (Euclidean) distance $\ell^\g$ of one of the $x_i$.
We may arrange that $n(\ell) \le \cb \ell^{(1-\g)(d-1)}$ as the 
following constructive sketch shows.  Take $x_1 = \ell \hat e_1$
and suppose $x_1,\dots,x_k$ have already been selected.
Choose $x_{k+1}\in \partial \B(0,\ell) 
\setminus (\cup_{i=1}^k \B(x_i, \ell^\g))$ if this latter set
is non-empty, and stop otherwise. The Eulcidean balls $\B(x_i,\ell^\g/2)$
cover disjoint patches of $\partial \B(0,\ell)$ with $(d-1)$-dimensional
area of order $\ell^{\g(d-1)}$.  Since $\partial \B(0,\ell)$ has
total area of order $\ell^{d-1}$, it follows that the process must stop
after order $\ell^{(1-\g)(d-1)}$ steps.  

Also, take $x'_i = 2\ell \hat e_1 - x_i$ so each $x'_i$ is on 
$\partial \B(2\ell \hat e_1, \ell)$ and 
every $x\in \partial \B(2\ell \hat e_1, \ell)$ is
within distance $\ell^\g$ of one of the $x'_i$. The $x'_i$ are
simply the $x_i$ radially reflected about $x_1 = \ell\hat e_1$ and
they bear the same spatial relation to each other as do the $x_i$.

We claim that for some constant $\cb$, for large $\ell$ we have
$$
T_{2\ell} \qge \min_{1\le i\le n(\ell)} T(0,x_i)\ +\ \min_{1\le j\le n(\ell)} T(2\ell\hat e_1,x'_j)\ -\ \cb\ell^{\g\a}
\text{\quad on $\F_\ell^c$}.\tag{\eqcd}
$$
To see this, let $r = (q_k)$ denote the path from $q(0)$ to $q(2\x)$
that realizes $T_{2\ell}$. 
Let $\tilde q = q_{k^*}$ denote the first
$q_k$ on $r$ not in  $\B(0,\ell)$ and put $q = q_{k^*-1}$. 
(Since $r$ ends with
$q(2\x)$ and, on $\F_\ell^c$, 
$|q(2\x) - 2\x| < \ell^\g < \ell$, such a $k^*$ exists; furthermore,
$k^*\not= 0$ since $r$ begins with $q(0)$ and $|q(0)| < \ell$ on $\F_\ell^c$.)
Similarly, let $q'$ denote the first $q_k$ on $r$ such that
$q'$ and all subsequent $q_k$'s on $r$ lie within $\B(2\x,\ell)$.
Then clearly 
$$
T_{2\ell} \qge T(0,q) + T(q',2\x).
$$

Now let $x = \overline{q\,\tilde q}\,\cap\,\partial \B(0,\ell)$; it
follows from (5.3) of Lemma 5.2 that,
for some $\cb$, on $\F_\ell^c$ we must have $|q-x| \le \cb\ell^\gamma$
for all large $\ell$. 
Picking $x_{i^*}$ so that $|x_{i^*} - x| \le \ell^\g$, we get that
$$
|q(x_{i^*}) - q| \qle |q(x_{i^*}) - x_{i^*}| + |x_{i^*} - x| + |x - q|
\qle (2+\cb)\ell^\g.
$$
It follows that $T(0,x_{i^*}) \le T(0,q) + (2+\cb)^\a\ell^{\a\g}$ and hence
$$
T(0,q) \qge \min_{1\le i\le n(\ell)} T(0,x_i) - (2+\cb)^\a\ell^{\a\g}.
$$
Similarly, 
$$
T(2\x,q') \qge \min_{1\le j\le n(\ell)} T(2\x,x'_j) - (2+\cb)^\a\ell^{\a\g},
$$
yielding (\eqcd)
for an appropriately larger $\cb$.
Since $x_1 = x'_1 = \x$, it follows that
$$
\min_{1\le i\le n(\ell)} T(0,x_i)\ +\ \min_{1\le j\le n(\ell)} 
T(2\ell\hat e_1,x'_j) \qle T_{2\ell} + 
\cb\ell^{\a\g} + T(0,\x)I_{\F_\ell} + T(2\x,\x)I_{\F_\ell}.
$$

Taking expectations and using the symmetry of our construction together
with the Cauchy-Schwarz  inequality yields
$$
2E[\min_i T(0,x_i)] \qle ET_{2\ell} + 
\cb\ell^{\a\g} + 2\sqrt{E[T_\ell^2]P[\F_\ell]}. \tag{\eqce}
$$
Now $E[T_\ell^2] = (ET_\ell)^2 + \var T_\ell$, where the second summand
is of order $\ell$ by Theorem 2.1 and the first term is of order $\ell^2$
by general subadditivity arguments (see (7) in [HoN1]). 

It follows from (5.2) of Lemma 5.2 
(for possibly different $\cs$ and $\cb$)
that $P[\F_\ell] \le \cb\exp(-\cs \ell^{\g d})$.
Hence 
$$
\sqrt{E[T_\ell^2]P[\F_\ell]} \qequals o(1) \qequals o(\ell^{\a\g}) \text{\ as $\ell\to\infty$}
$$
and 
$$
\aligned
ET_{2\ell}&\qge 2E[\min_{1\le i\le n(\ell)} T(0,x_i)]  - \cb\ell^{\a\g}\\
          &\qequals 2ET_\ell\ -\ 2E[\max_{1\le i\le n(\ell)} 
            \Big(E[T(0,x_i)] - T(0,x_i)\Big)] \ -\ \cb\ell^{\a\g}.\\
\endaligned
$$
The equality above uses that $E[T(0,x_i)] = E[T(0,x_1)] = ET_\ell$.
Since $\a\g<1/2$, Lemma 4.1 will be proved if we establish that 
$$
E[\max_{1\le i\le n(\ell)}\Big(E[T(0,x_i)] - T(0,x_i)\Big)] \qle 
\cb\sqrt{\ell}(\log\ell)^{1/\kappa_1}. \tag{\eqcf}
$$
To conclude the proof of Lemma 4.1, take 
$Y_i^{(\ell)} = T(0,x_i)/\sqrt{\ell}$
in Lemma 4.3 below and note that 
hypotheses are satisfied with $a = \frac 1 2 + \epsilon$,
$\tilde a = (1-\gamma)(d-1) + \epsilon$, $b = \kappa_1$, $\tilde b = \kappa_2$, and $\cs$ and $\cb$ as in Theorem 2.1.

\proclaim{Lemma 4.3} For $\ell\ge \ell_0>1$,
let $Y^{(\ell)}_i$ for $1\le i\le n(\ell)$
be non-negative
random variables on a common probability space such that, for some
$a, \tilde a, b, \tilde b, \cs, \cb \in (0,\infty)$, 
$$
E[Y^{(\ell)}_i] \qle \ell^a 
\text{\quad and\quad} n(\ell)\qle\ell^{\tilde a}, \tag{\eqcg}
$$ 
and, 
$$
\text{$P(|Y^{(\ell)}_i - E[Y^{(\ell)}_i]| > x) \qle 
\cb\exp(-\cs x^b)$\quad for $x\le \cs\ell^{\tilde b}$} \tag{\eqch}
$$
Then, for some $C_2 = C_2(\ell_0, a,\tilde a,b, \tilde b,\cs,\cb)$, 
$$
E[\max_{1\le i\le n(\ell)}(E[Y^{(\ell)}_i] - Y^{(\ell)}_i)] 
\qle C_2(\log \ell)^{1/b}\quad\text{for all $\ell\ge\ell_0$}. \tag{\eqci}
$$
\endproclaim

\demo{Proof}  Let $M^{(\ell)}$ denote 
$\max_{1\le i\le n(\ell)}(E[Y^{(\ell)}_i] - Y^{(\ell)}_i)$ and
put $f(\ell) = \hat C(\log \ell)^{1/b}$ where we take
$\hat C$ so that $\cs\hat C^b = a+\tilde a$.
Note that $M^{(\ell)} \le \ell^a$ since the $Y^{(\ell)}_i$ are non-negative, so
$$
M^{(\ell)} \qle \cases    f(\ell) &\text{if $Y^{(\ell)}_i - E[Y^{(\ell)}_i] \ge -f(\ell)$ for all $i\le n(\ell)$}\\
                       \ell^a&\text{otherwise.}\\
              \endcases
$$
For large $\ell$, 
$f(\ell) \le \cs\ell^{\tilde b}$ and we have
$$
\align
E[M^{(\ell)}] &\qle f(\ell)\ +\ \ell^a\sum_{i=1}^{n(\ell)} 
P\Big(Y^{(\ell)}_i - E[Y^{(\ell)}_i] \le -f(\ell)\Big)\\
&\qle f(\ell)\ +\ \ell^{a+\tilde a}\cb\exp(-\cs f(\ell)^b)\\
&\qequals f(\ell)\ +\ \cb \qle C_2 (\log \ell)^{1/b},\\
\endalign
$$
where the equality follows from our choice of $\hat C$ and the final 
inequality holds for an appropriate $C_2$.  The second inequality above
holds only
for large $\ell$, but since $EM^{(\ell)} \le \ell^a$ we can ensure that 
$EM^{(\ell)} \le C_2(\log \ell)^{1/b}$ for all $\ell\ge \ell_0$ by
making $C_2$ larger if necessary.\qed
\enddemo

\enddemo

\head{5. Technical Lemmas}\endhead

Throughout this section, $\phi$ is any cost function of the form
$$
\phi(t) \qequals \cases 
                      t^\alpha &\text{if $t \le h$}\\
                      h^\alpha + \alpha h^{\alpha-1}(t-h)&\text{otherwise},
                      \endcases
$$
with $\alpha > 1$ and $h>0$.
Recall our notation that, 
for any cost function $\phi$ of this form and $a,b\in\rd$, 
$$
\W_\phi(a,b) \qequals \{c\in\rd: \phi(|a-c|) + 
\phi(|c-b|) \le \phi(|a-b|)\} \tag5.1
$$
and that $\W(a,b) = \W_{\phi_\infty}(a,b)$, where $\phi_\infty(t) = t^\a$.  
We provide below in Lemma 5.1 some elementary geometric properties
of these regions. 

\proclaim{Lemma 5.1} The region $\W_\phi(0,\x)$ is closed and convex,
contains $\frac1 2 \x$ in its interior, and is invariant with respect to rotations about the first coordinate axis.   
Also, $\W_\phi(a,b)$ is the set $\W_\phi(0,|a-b| \hat e_1)$
rigidly moved so that $0$ is moved to $a$ and $|a-b|\hat e_1$ is moved to $b$. (By the rotational
invariance of $\W_\phi(0,|a-b|\hat e_1)$ about the first coordinate axis, 
any such rigid motion will do.)
In the case $\phi = \phi_\infty$, 
$\W(0,\x) = \ell \W(0,\hat e_1)$ and $\ell' < \ell$ implies that
$\W(0,\ell' \hat e_1) \subset \W(0,\x)$. 
\endproclaim

\demo{Proof} Much of this Lemma is self-evident. We prove only
the convexity claim and the statements about the case $\phi = \phi_\infty$.
The convexity of
${\Cal W}_\phi (0,\x)$ follows 
from the facts that $\phi$ is
convex and increasing as follows.  For $c,c'\in \W_\phi(0,\x)$, and $\l\in[0,1]$:
$$
\align
\phi(\ell) &\qge \l(\phi(|c|) + \phi(|c-\x|))\ +\ (1-\l)(\phi(|c'|) + \phi(|c'-\x|))\\
&\qge \phi(\l|c|+(1-\l)|c'|)\ +\ \phi(\l|c-\x|+(1-\l)|c'-\x|)\\
&\qge \phi(|\l c+(1-\l)c'|)\ +\ \phi(|\l(c-\x)+(1-\l)(c'-\x)|)\\
&\qequals \phi(|\l c+(1-\l)c'|)\ +\ \phi(|\l c+(1-\l)c' - \x)|),\\
\endalign
$$
so also $\l c+(1-\l)c'\in{\Cal W}_\phi (0,\x)$.
That 
${\Cal W}(0,\x) = \ell {\Cal W}(0,\hat e_1)$ follows from
the (degree $\a$) homogeneity  of $\phi_\infty$. If $\ell'<\ell$,
${\Cal W}(0,\ell'\hat e_1)
=\ell' {\Cal W}(0,\hat e_1)\subset\ell {\Cal W}(0,\hat e_1)
={\Cal W}(0,\x)$, where the containment follows 
since $0$ is in the convex ${\Cal W}(0,\hat e_1)$.\qed
\enddemo

\proclaim{Lemma 5.2} For $\g\in(0,1)$, let
$
A_{\g,\ell} \equiv \{\text{$\exists$ $a\in\rd$ with $|a|\le 2\ell$ and $|a-q(a)|\ge\ell^\g$}\}.
$
Then, for some $\cs$ and $\cb$:
$$
P[A_{\g,\ell}]\le \cb\exp(-\cs\ell^{\g d}),\tag5.2
$$
and furthermore, for large $\ell$, on $A_{\g,\ell}^c$,
$$
\sup\{|a-b|:|a|\le\ell,\ b\in\rd,\ \W(a,b)\cap Q = \emptyset\} \qle \cb\ell^\g.\tag5.3
$$
\endproclaim

{\bf Remark.}  If $\Gamma\equiv\sup\{|a| : 
a\in\rd,\ \W(0,a)\cap Q = \emptyset\}$, then (for large $\ell$) 
$\Gamma \le \cb\ell^\g$ on $A^c_{\g,\ell}$.
By the substitution $x = \cb\ell^\g$, it follows that
$$
P[\Gamma>x]\qle \cb\exp(-\cs x^d)\tag5.4
$$ 
(for possibly different $\cs$ and $\cb$).
Also, on $A^c_{\g,\ell}$, if $(q,q')$ is any geodesic segment
with $|q|\le\ell$ (or $|q'|\le\ell$), then $|q-q'|\le \cb \ell^\g$.
It follows (for possibly different $\cs$ and $\cb$) that 
$$
P[\text{$\exists$ geodesic segment $(q,q')$ with $|q|\le\ell$ or $|q'|\le\ell$
 and $|q-q'|>\ell^\g$}] \qle \cb\exp(-\cs\ell^{\g d}).\tag5.5
$$
While Lemma 5.2 gives (5.4) and (5.5) for large $x$ and $\ell$, respectively, this restriction is removed by increasing $\cb$.

\demo{Proof of Lemma 5.2} For large $\ell$, we have that
$$
A_{\g,\ell} \ \subset\ \{\text{$\exists$ $a\in\zd$ with $|a|\le 2\ell$ and $|a-q(a)|\ge\ell^\g/2$}\}.
$$
This larger event has probability bounded by 
$\cb\ell^d\exp(-\cs\ell^{\g d})$, which, for smaller $\cs$
is bounded (for large $\ell$) by $\cb\exp(-\cs\ell^{\g d})$.
By increasing $\cb$ if necessary, (5.2) will hold for all $\ell$.

To get (5.3), we take $\cb$ large enough so that 
$$
{\Cal B}(\frac1 2\hat e_1, \cb^{-1})\ \subset\ \W(0,\hat e_1).\tag5.6
$$
Suppose $\ell$ is large enough so that $\ell>\cb\ell^\g$ and so that,
for a configuration $Q$, we can find
$a,b\in\rd$ satisfying: $|a|\le \ell$, $|a-b| > \cb\ell^\g$, with
$\W(a,b)$ devoid of particles from $Q$. If $|b| < 2\ell$,
but $\tilde b = b$; otherwise put $\tilde b = \overline{a\,b}\,\cap\,\partial{\Cal B}(0,2\ell)$.
Then, since $|a-\tilde b|\ge\ell$ and 
$\W(a,\tilde b)\subset \W(a,b)$ (by Lemma 5.1), $a$ and $\tilde b$ satisfy:
$|(a+\tilde b)/2| < 2\ell$, $|a-\tilde b|>\cb\ell^\g$, with 
$ \W(a,\tilde b)$ devoid of Poisson particles.
Since, using (5.6), 
${\Cal B}((a+\tilde b)/2, \cb^{-1}|a-\tilde b|)\subset \W(a,\tilde b)$,
it follows that 
$$
\bigg|q\big(\frac{a+\tilde b}{2}\big) - \frac{a+\tilde b}{2}\bigg|
\qge \cb^{-1}|a-\tilde b| \qg \ell^\g,
$$
i.e., the configuration $Q$ belongs to $A_{\g,\ell}$.\qed

\enddemo

\proclaim{Lemma 5.3}  For any $a,b,c\in\rd$ we have:
$$
\phi^2(|a-c|) \qle 2^{2\a}(\phi^2(|a-b|) + \phi^2(|b-c|))\tag{\eqbb}
$$
and
$$
\phi(|a-c|) - \phi(|a-b|) - \phi(|b-c|) \qle 2^\a h^\a.\tag{\eqbc}
$$
\endproclaim

\demo{Proof}  First we prove (\eqbb).  If $t \le 2h$ then
$$
\frac{\phi(t)}{\phi(t/2)} \qequals \frac{\phi(t)}{(t/2)^\a} \qle \frac{t^\a}{(t/2)^\a} \qequals 2^\a.
$$
If $t>2h$, put $t = (1+y)2h$ where $y>0$. Then
$$
\frac{\phi(t)}{\phi(t/2)} \qequals \frac{1+\a(1+2y)}{1+\a y} \qequals 1 + \a\frac{1+y}{1+\a y}
\qle 1+\a\qle 2^\a,
$$
with the latter two inequalities holding since $\a>1$.
Thus $\phi (t) \le 2^\a \phi(t/2)$ for all $t \ge 0$.
Now suppose, without loss of generality, that $|a-b|\le|b-c|$ so $|b-c|\ge\frac1 2 |a-c|$ and
$$
\phi(|a-c|)\qle 2^\a \phi(\frac1 2 |a-c|) \qle 2^\a\phi(|b-c|)
$$
giving that
$$
\phi^2(|a-c|) \qle 2^{2\a}\phi^2(|b-c|) 
\qle 2^{2\a}(\phi^2(|a-b|)+\phi^2(|b-c|)),
$$
and verifying (\eqbb). To establish (\eqbc), 
we first show by examining cases that, for $x,y\ge0$, 
$$
\phi(x+y) -\phi(x) -\phi(y) \le 2^\a h^\a.\tag{\eqbd}
$$
This clearly holds if $x+y\le 2h$. If $x+y>2h$ 
with $x>h$ and $y\le h$, then 
$$
\phi(x+y)-\phi(x) -\phi(y)\qle \phi(x+y)-\phi(x) \qequals \a h^{\a-1}y \qle \a h^\a \qle 2^\a h^\a.
$$
A symmetric argument works for $x+y>2h$ with $x\le h$ and $y>h$. 
Finally, if $x>h$ and $y>h$, then
$$
\phi(x+y)-\phi(x) -\phi(y) \qequals (\a-1) h^\a \qle 2^\a h^\a.
$$
To complete the proof of (\eqbc), 
let $b'$ be the orthogonal projection of $b$ onto the line passing through $a$ and $c$.
Then the left side of (\eqbc) is dominated by 
$\phi(|a-c|) - \phi(|a-b'|) - \phi(|b'-c|)$.  If $b'\not\in \overline{ac}$, this
quantity is negative.  If $b'\in\overline{ac}$, then (\eqbd) yields (\eqbc).\qed

\enddemo

\proclaim{Lemma 5.4}  For any $E>0$ and $a,b\in\rd$, let ${\Cal H}_E(a,b)$ denote the set
$$
\align
{\Cal H}_E(a,b) \qequals \{c\in\rd:&\text{ $\exists$ a point $p$ on the line segment connecting}\\
                              &\text{${\scriptstyle\frac3 4} a + {\scriptstyle\frac1 4} b$ 
                                 and ${\scriptstyle\frac1 4} a + {\scriptstyle\frac3 4} b$
                                     such that $|c-p|\le E$}\},\\
\endalign
$$
and define $\W_\phi(a,b)$ as in (5.1).
Then for any $E>0$, there is an $h_0>0$ such that 
${\Cal H}_E(a,b)\subset  \W_\phi(a,b)$
whenever $|a-b|>h_0$ and $h>h_0$.
\endproclaim

\demo{Proof}
Clearly it suffices to prove this for $a=0$ and $b=\x$ where $\ell>0$.  Let $c$ be any
point whose $\hat e_1$ coordinate is $\ell/2$ and put 
$u = |c - (\ell/2)\hat e_1|$.
First, by examining cases, we calculate how large $u$ may be while
keeping $c$ inside $ \W_\phi(0, \x)$.
Since $|c| = |c-\x|$, to have $c\in \W_\phi(0,\x)$ we need $2\phi(|c|) \le \phi(\ell)$ for 
which it is sufficient to have:
$$
2\phi\big({\frac \ell 2} + u\big) \qle \phi(\ell).\tag{\eqbe}
$$
If $\ell<h$, to have (\eqbe), it 
suffices to have $2(\frac \ell 2 + u)^\a \le \ell^\a$
or 
$$
u\le (2^{-1/\a} - 2^{-1}) \ell.\tag{\eqbf}
$$
On the other hand, if $\ell > 2h$, (\eqbe) will obtain provided
$$
2(h^\a + \a h^{\a-1}({\frac \ell 2} + u - h)) \qle h^\a + \a h^{\a-1}(\ell-h),
$$
which reduces to:
$$
u\qle \frac{\a-1}{2\a} h.\tag{\eqbg}
$$
Finally, if $h\le \ell\le2h$, it suffices to have
$$
2\big({\frac \ell 2} + u\big)^\a \qle h^\a + \a h^{\a-1}(\ell-h),
$$
or, equivalently,
$$
u\qle \bigg[2^{-1/\a} (1 + \a({\frac \ell h}-1))^{1/\a} - {\frac1 2 \frac \ell h}\bigg]  h.
$$
One verifies by calculus that the quantity in brackets, viewed as a function of $\ell$,
is increasing on the interval $[h,\frac{\a+1}{\a} h]$ and decreasing on
$[\frac{\a+1}{\a} h, 2h]$. It therefore suffices for the case
$h\le \ell\le 2h$ to have 
$$
u\qle \min(2^{-1/\a} - 2^{-1}, 2^{-1/\a}(1+\a)^{1/\a} - 1) h. \tag{\eqbh}
$$
(Note that this minimum is strictly greater than 0 since $\a>1$.)
Using (\eqbf), (\eqbg), and (\eqbh), we see that 
to ensure that $c\in  \W_\phi(0,\x)$
it suffices to have  
$$
\aligned
u&\qle C \min(\ell, h),\text{\ where}\\
C &\qequals \min({\frac{\a-1}{2\a}},\ 2^{-1/\a}-2^{-1},\ 2^{-1/\a}(1+\a)^{1/\a} - 1).\\
\endaligned
$$
That is,
$$
{\Cal U}_E(\ell) \qequals \{c\in\rd : \text{$c$'s $\hat e_1$ coordinate is $\ell/2$, 
         and $|c-(\ell/2)\hat e_1| \le E$}\}
$$
satisfies ${\Cal U}_E(\ell) \subset \W_\phi(0,\x)$ if $E\le C\min(\ell,h)$.
It follows from the convexity of $ \W_\phi(0,\x)$ that the suspension of ${\Cal U}_E(\ell)$
defined by
$$
{\Cal S}_E(\ell) \qequals \{\rho\, {\Cal U}_E(\ell) : 0\le\rho\le 1\} \cup
                             \{\rho\, {\Cal U}_E(\ell) + (1-\rho)\x : 0\le\rho\le1\},
$$
also satisfies 
$$
{\Cal S}_E(\ell) \ \subset \W_\phi(0,\x)\text{\ for $E\le C\min(\ell,h)$}. \tag{\eqbi}
$$
Elementary geometric arguments show that 
${\Cal H}_E(\ell)\subset{\Cal S}_{4E}(\ell)$ if $\ell\ge8E$. 
It follows from this and (\eqbi) that
$$
{\Cal H}_E(\ell) \ \subset\ \W_\phi(0,\x)\text{\ for\ $\ell\ge8E$ and $(C/4)\min(\ell,h)\ge E$},
$$
proving the lemma for $h_0=\max(8E, 4E/C)$.\qed
\enddemo

The next purely geometric lemma (proved in [Ho]) 
states, roughly speaking that
if $(q_0,\dots,q_n)$ is a minimizing path with respect to
the cost function $\phi$ and a segment
$L = \overline{q_i\,q_{i+1}}$ passes near a segment 
$L' = \overline{q_{i'}\,q_{i'+1}}$ where $i< i'$,
then this must happen near the end of $L$ and the beginning of
$L'$. Specifically:

\proclaim{Lemma 5.5 (No Doubling Back Proposition [Ho])} Under the above
arrangement, if $a\in L$ and 
$b\in L'$,  then $|q_{i+1}-a|\le 16|a-b|$ and
$|q_{i'}-b|\le16|a-b|$.  Also, therefore, $|q_{i+1}- q_{i'}| \le 33|a-b|$.
\endproclaim

The following lemma is a modification of Theorem 3 of [Ke2].

\proclaim{Lemma 5.6}  Let $(M_k:k\ge0)$, $M_0\equiv0$, be a 
martingale with respect to the filtration $\field_k\uparrow \field$.
Put $\D_k = M_k - M_{k-1}$ and suppose $(U_k:k\ge1)$ 
is a sequence of $\field$-measurable positive random 
variables satisfying $E[\D_k^2|\field_{k-1}] \le E[U_k|\field_{k-1}]$.  
With $S=\sum_{k=1}^\infty U_k$, suppose further
that for finite 
constants $C'_1>0$, $0<\g\le1$, $c\ge1$, and $x_0\ge c^2$ we have
$|\D_k|\le c$ and
$$
P[S > x] \qle C'_1 \exp(-x^\g),\text{\ when $x\ge x_0$}.\tag\eqbj
$$
Then $\lim_{k\to\infty}M_k = M$ exists and is finite almost surely and
there are constants 
(not depending on $c$ and $x_0$) 
$C_2 = C_2(C'_1,\g)<\infty$ and $C_3 = C_3(\g)>0$ such that
$$
P[|M| \ge x \sqrt{x_0}\,] \qle C_2 \exp(-C_3 x)\text{\ when $x\le x_0^\g$}.
$$
\endproclaim

\demo{Proof} The proof of this lemma largely parallels the
proof of Theorem 3 of [Ke2].

Throughout the proof, $C_2(C'_1,\g)$ will denote a constant whose value
depends only on $C'_1$ and $\g$.  As the proof progresses, $C_2$ will be made possibly larger 
several times, each occurrence of which is indicated by a ``$+$'' superscript: $C_2^+(C'_1,\g)$.
Similarly, $C_3(\g)$ will be made possibly smaller 
when indicated by a ``$-$'' superscript.

Following Kesten, put 
$$\aligned
A &\qequals \sum_{k=1}^\infty E[\D_k^2|\field_{k-1}],\cr
\nu &\qequals \inf\Big\{\ell:\sum_{k=\ell+1}^\infty E[U_k|\field_\ell] > z\Big\}
\text{\ (where $\inf\emptyset = \infty$), and}\cr
\tilde A &\qequals \sum_{k=1}^\nu E[\D_k^2|\field_{k-1}].
\endaligned
$$
Here $z > 0$ is arbitrary, but a specific choice will be made later.
Then it follows exactly as in Kesten's Step 2 that
$$
P[A\ge y] \qle P[\nu<\infty]\ +\ P[\tilde A \ge y]\tag\eqbk
$$
and that, for any positive integer $r$, 
$E[\tilde A^r] \qle r!z^{r-1} ES$.
Next, we estimate
$$
ES \qequals \int_0^\infty P[S > s]\,ds
\qle x_0 + C'_1\int_0^\infty \exp(-s^\g)\,ds
\qequals x_0 + C_2(C'_1,\g),
$$
so
$E[\tilde A^r] \qle r!z^{r-1} (x_0 + C_2(C'_1,\g))$.
Also, as in Kesten's (5.8), by taking $r=\lfloor y/z\rfloor$ where $y\ge z$, we get
$$
\align
P[\tilde A \ge y] &\qle C'\cdot (x_0 + C_2(C'_1, \g))\frac1 z \exp(-\frac y {2z})\cr
&\qle C_2^+(C'_1, \g) \exp(-\frac y {2z})\tag\eqbl\cr
\endalign
$$
with the second inequality holding for $y\ge z\ge x_0$ since 
also $x_0\ge1$. ($C'$ comes
from Stirling's formula and the fact that 
$\sqrt{y/z} \le \text{constant} \cdot \exp(y/2z)$.)

Next, as in Kesten's Step 3, we estimate $P[\nu < \infty]$.
Let $S_m = \sum_{k=1}^m U_k$ and $S_{m,\ell} = E[S_m | \field_\ell]$.
If $g(s) = \exp(\frac1 2 s^\g)$ then $g'(s) = \frac1 2 \g s^{\g-1}\exp(\frac1 2 s^\g)>0$
for $s>0$ and $g''(s) = \frac1 2 \g s^{\g-2}\exp(\frac1 2 s^\g)(\frac1 2 \g s^\g +\g-1)>0$
when $s^\g>2(1-\g)/\g = 2\b$.  
Hence 
$\tilde g(s) =(e^\b \vee \exp(\frac1 2 s^\g))$ is convex
giving that $\big(\tilde g(S_{m,\ell}): \ell\ge0\big)$ is a submartingale.
Also, for $z\ge z(\g) = (2\b)^{1/\g}$, $\tilde g(s) > \exp(\frac1 2 z^\g)$ if and only if
$s > z$. So, for $z \ge z(\g)$:
$$\align
P[\nu < \infty] &\qle \lim_{m\to\infty} \lim_{n\to\infty} P[\max_{\ell\le n} S_{m,\ell} > z]\cr
&\qequals \lim_{m\to\infty} \lim_{n\to\infty}P\Big[\max_{\ell\le n} \big\{\tilde g(S_{m,\ell})\big\} > \exp(\frac1 2 z^\g)\Big]\cr
&\qle \limsup_{m\to\infty} \limsup_{n\to\infty}\exp(-\frac1 2 z^\g)\, E[\tilde g(S_{m,n})]\text{\ (by Doob's inequality)}\cr
&\qle \limsup_{m\to\infty}\limsup_{n\to\infty}\exp(-\frac1 2 z^\g)\, E[E\{\tilde g(S_m) | \field_n\}]\text{\ (by Jensen's inequality)}\cr
&\qle \limsup_{m\to\infty}\exp(-\frac1 2 z^\g)\, \big(e^\b + E[g(S_m)]\big)\cr
&\qequals \exp(-\frac1 2 z^\g)\, \big(e^\b + E[g(S)]\big).\cr
\endalign
$$
Now 
$$
\align
E[g(S)] &\qle g(x_0)P[S\le x_0]\ -\ \int_{x_0}^\infty g(s)\,dP[S> s]\cr
&\qequals g(x_0)\ +\  \int_{x_0}^\infty g'(s) P[S>s]\,ds,\cr
&\qle \exp(\frac1 2 x_0^\g)\ +\ C'_1\frac{\g}{2}\int_{x_0}^\infty \exp(-\frac1 2 s^\g)\,ds\text{\ \ (since $s^{\g-1}\le1$ on $[x_0,\infty)$)}\cr
&\qle \exp(\frac1 2 x_0^\g) \ +\ C_2^+(C'_1,\g)
\text{\ \ (by replacing $\int_{x_0}^\infty$ with $\int_0^\infty$).}\cr
\endalign
$$
Hence, for $z > z(\g)$,
$$
\align
P[\nu < \infty] &\qle \exp(-\frac1 2 z^\g) \big(C_2^+(C'_1,\g) +  \exp(\frac1 2 x_0^\g)\big)\cr
&\qle C_2^+(C'_1,\g) \exp\big(-\frac1 2 (z^\g - x_0^\g)\big).\tag\eqbm\cr
\endalign
$$
Following Kesten again by letting $y\to\infty$ and 
then $z\to\infty$, (\eqbl), (\eqbm), and (\eqbk) give
that $P[A = \infty] = 0$.  But 
$\lim_{k\to\infty} M_k = M$ exists and is finite almost
surely on $\{A<\infty\}$ (See, e.g., Theorem 4.8 of [D].)

Next, as in Step 1 of Kesten and pp. 154-155 of [Ne]
(this is where the boundedness of the martingale differences is used), 
for $y\ge cx>0$:
$$
P[M\ge x] \qle P[A\ge y]\ +\ \exp(-\frac{x^2}{2ey}).\tag\eqbn
$$
Combining (\eqbk), (\eqbl), (\eqbm), and (\eqbn), we get that
$$
P[M\ge x] \qle C_2^+(C'_1,\g) \Big[\exp(-\frac{z^\g - x_0^\g} 2) + \exp(-\frac y {2z}) + \exp(-\frac{x^2}{2ey})\Big]
$$
whenever 
$$
y \ge cx,\ y\ge z\ge x_0, \text{ and } z \ge z(\g).\tag\eqbo
$$

Now, like in Kesten's Step 4, take $z = (x_0^\g + x^a)^{1/\g}$ where 
$a = 2\g/(1+2\g)$, and $y = xz^{1/2}$. 
Then $2z^{1/2} \le 2^{1/\g} (x_0^{1/2} + x^{a/(2\g)})$ so, 
with $C_3(\g) = 2^{-1/\g}/e$, 
$$
\frac y {2z} = \frac x {2 z^{1/2}} \ge C_3(\g)\frac x {x_0^{1/2} + x^{a/(2\g)}}
\text{\ and }
\frac {x^2} {2ey} = \frac x {2 e z^{1/2}} \ge C_3(\g)\frac x {x_0^{1/2} + x^{a/(2\g)}}.
$$
Also, since $a = 1 - a/(2\g)$ and $C_3(\g) < 1/2$, 
$$
(z^\g - x_0^\g)/2 \qequals  x^a/2 \qequals 
\frac {x/2} {x^{a/(2\g)}} \qge C_3(\g)\frac x {x_0^{1/2} + x^{a/(2\g)}}.
$$
Presently we verify that for some constant $C_4(\g)$, (\eqbo) holds provided $x \ge C_4(\g)\sqrt{x_0}$.
The relation $y\ge cx$ is equivalent to $z\ge c^2$; but $z\ge x_0 \ge c^2$, giving two
inequalities in (\eqbo).  To get $z\ge z(\g) = (2\b)^{1/\g}$, it suffices to have
$x\ge (2\b)^{1/a}$ which, since $x_0\ge c^2\ge1$, will hold if $x\ge (2\b)^{1/a}\sqrt{x_0}$.
Finally, $y\ge z$ is equivalent to $x^{2\g} \ge x_0^\g + x^a$ which will hold provided
$$
\frac1 2 x^{2\g} \ge x_0^\g\text{\ and } \frac1 2 x^{2\g}\ge x^a,
$$
or, equivalently, when
$$
x \ge 2^{1/2\g} \sqrt{x_0}\text{\ and } x \ge 2^{1/(2\g-a)}.\tag\eqbp
$$
Since $1/(2\g-a) = (1+2\g)/4\g^2 \ge 1/2\g$ and $x_0\ge1$, both conditions in (\eqbp) will hold
provided $x \ge 2^{(1+2\g)/4\g^2}\sqrt{x_0}$.  It therefore suffices to take
$C_4(\g) = \max((2\b)^{1/a}, 2^{(1+2\g)/(4\g^2)})$.

Letting $d = d(\g) = 2\g + 1$, we get
$$
P[M \ge x] \qle C_2^+(C'_1,\g) \exp\bigg[-C_3(\g)\frac x {x_0^{1/2} + x^{1/d}}\bigg]
$$
whenever $x \ge C_4(\g) \sqrt{x_0}$.
Now, for $C_4(\g)x_0^{1/2} \le \tilde x \le x_0^{d/2}$ we also have $\tilde x^{1/d} \le x_0^{1/2}$, so
$$
P[M \ge \tilde x] \qle C_2(C'_1,\g) \exp\bigg[-C_3^-(\g) \frac {\tilde x} {\sqrt{x_0}}\bigg].
$$
Substituting $\tilde x = x\sqrt{x_0}$, we get that, for $C_4(\g)\le x\le x^\g_0$, 
$$
P[M \ge x\sqrt{x_0}\,] \qle  C_2(C'_1,\g) \exp[-C_3(\g) x].
$$
But for $x < C_4(\g)$, the exponential is bounded away from zero by $\exp(-C_3(\g)C_4(\g))$.
Hence,
$$
P[M \ge x\sqrt{x_0}\,] \qle  C_2^+(C'_1,\g) \exp[-C_3(\g) x] \text{\ for $x \le x^\g_0$}.
$$
The lemma follows by a further 
application of this to the martingale $(-M_k: k\ge 0)$.
\qed
\enddemo

\proclaim{Acknowledgment} We thank the 
anonymous referee for an impressively complete
report, including very useful suggestions that improved the presentation
of the paper.
\endproclaim


\Refs


\refstyle{A}
\widestnumber\key{CGGKXX}

\ref \key A
     \by Aizenman, M.
     \paper The geometry of critical percolation and conformal invariance
     \inbook The 19th IUPAP International Conference on Statistical
     Physics
     \ed H. Bai-lin
     \publ World Scientific
     \yr 1996
     \publaddr Singapore
     \pages 104-120
\endref\vskip.08truein

\ref \key AB
     \by Aizenman, M. and Burchard, A.
     \paper H\"{o}lder regularity and dimension bounds for random curves
     \jour Duke Math. J.
     \yr 1999
     \vol 99
     \pages 419-453
\endref\vskip.08truein

\ref \key ABNW
     \by Aizenman, M., Burchard, A., Newman, C.M. and Wilson, D.
     \paper Scaling limits for minimal and random spanning trees in
     two dimensions
     \jour Random Struct. Alg.
     \yr 1999
     \vol 15
     \pages 319-367
\endref\vskip.08truein

\ref \key AS
     \by Aldous, D. and Steele, J.M.
     \paper Asymptotics for Euclidean minimal spanning trees on random points
     \jour Probab. Th. Rel. Fields
     \yr 1992
     \vol 92
     \pages 247-258
\endref\vskip.08truein

\ref \key Al1
     \by Alexander, K.S.
     \paper A note on some rates of convergence in first-passage percolation
     \jour Ann. Appl. Probab.
     \vol 3
     \yr 1993
     \pages 81-90
\endref\vskip.08truein

\ref \key Al2
     \by Alexander, K.S.
     \paper Percolation and minimal spanning trees in infinite graphs
     \jour Ann. Probab.
     \vol 23
     \yr 1995
     \pages 87-104
\endref\vskip.08truein

\ref \key Al3
     \by Alexander, K.S. 
     \paper Approximations of subadditive functions and convergence\break rates in limiting-shape
            results
     \jour Ann. Probab.
     \vol 25
     \yr 1997
     \pages 30-55
\endref\vskip.08truein



\ref \key AlM
     \by Alexander, K.S. and Molchanov, S.A.
     \paper Percolation of level sets for two-dimensional random fields with lattice symmetry
     \jour J. Statist. Phys.
     \vol 77
     \yr 1994
     \pages 627-643
\endref\vskip.08truein

     

\ref \key BDJ
     \by Baik, J. Deift, P. and Johansson, K.
     \paper On the distribution of the longest increasing subsequence in
     a random permutation
     \jour J. Amer. Math. Soc. 
     \vol 12
     \yr 1999
     \pages 1119-1178
\endref\vskip.08truein


\ref \key BK
     \by van den Berg, J. and Kesten, H.
     \paper Inequalities with applications to percolation and reliability
     \jour J. Appl. Probab.
     \vol 22
     \yr 1985
     \pages 556-569
\endref\vskip.08truein

\ref \key Bo
     \by Boivin, D.
     \paper First-passage percolation:  The stationary case
     \jour Probab. Th. Rel. Fields
     \vol 86
     \yr 1990
     \pages 491-499
\endref\vskip.08truein



\ref \key CD
     \by Cox, J.T. and Durrett, R.
     \paper Some limit theorems for percolation processes with
            necessary and sufficient conditions
     \jour Ann. Probab.
     \yr 1981
     \vol 9
     \pages 583-603   
\endref\vskip.08truein


\ref \key CGGK
     \by Cox, J.T., Gandolfi, A., Griffin, P.S. and Kesten, H.
     \paper Greedy lattice animals I:  upper bounds
     \jour Ann. App. Probab.
     \yr 1993
     \vol 3
     \pages 1151-1169
\endref\vskip.08truein







\ref \key D
     \by Durrett, R.
     \book Probability: Theory and Examples
     \publ Wadsworth
     \publaddr Pacific Grove
     \yr 1991
\endref\vskip.08truein

\ref \key GK
     \by Gandolfi, A. and Kesten, H.
     \paper Greedy lattice animals II:  linear growth
     \jour Ann. App. Probab.
     \yr 1994
     \vol 4
     \pages 76-107
\endref\vskip.08truein

\ref \key Gr
     \by Grimmett, G.
     \book Percolation
     \publ Springer
     \publaddr Berlin-Heidelberg-New York
     \yr 1989
\endref\vskip.08truein

\ref \key HW
     \by Hammersley, J.M. and Welsh, D.J.A.
     \pages 61-110
     \paper First-passage percolation, subadditive processes, stochastic
       	   networks and generalized renewal theory
     \yr 1965
     \inbook Bernoulli, Bayes, Laplace Anniversary Volume
     \eds J. Neyman and L. LeCam
     \publ Springer
     \publaddr Berlin-Heidelberg-New York
\endref\vskip.08truein

\ref \key Ho
     \by Howard, C.D.
     \paper Good paths don't double back
     \jour Am. Math. Mon.
     \yr 1998
     \vol 105
     \pages 354-357
\endref\vskip.08truein

\ref \key HoN1
     \by Howard, C.D., and Newman, C.M.
     \paper Euclidean models of first-passage percolation
     \jour Probab. Th. Rel. Fields
     \yr 1997
     \vol 108
     \pages 153-170
\endref\vskip.08truein

\ref \key HoN2
     \by Howard, C.D., and Newman, C.M.
     \paper From greedy lattice animals to Euclidean first-passage percolation
     \inbook Perplexing Problems in Probability
     \eds M.~Bramson and R.~Durrett 
     \publ Birkh\"auser 
     \publaddr Boston-Basel-Berlin
     \yr 1999
     \pages 107-119
\endref\vskip.08truein

%

\ref \key HuH
     \by Huse, D.A. and Henley, C.L.
     \paper Pinning and roughening of domain walls in Ising systems due to random impurities
     \jour Phys. Rev. Lett.
     \vol 54
     \yr 1985
     \pages 2708-2711
\endref\vskip.08truein

\ref \key HuHF
     \by Huse, D.A., Henley, C.L. and Fisher, D.S.
     \jour Phys. Rev. Lett.
     \vol 55
     \yr 1985
     \pages 2924-2924
\endref\vskip.08truein

\ref \key J
     \by Johansson, K.
     \paper Transversal fluctuations for increasing subsequences on the plane
     \jour Probab. Theory Relat. Fields
     \vol 116
     \yr 2000
     \pages 445-456
\endref\vskip.08truein

\ref \key K
     \by Kardar, M.
     \paper Roughening by impurities at finite temperatures
     \jour Phys. Rev. Lett.
     \vol 55
     \yr 1985
     \pages 2923-2923
\endref\vskip.08truein

\ref \key KPZ
     \by Kardar, M., Parisi, G. and Zhang, Y.-C.
     \paper Dynamic scaling of growing interfaces
     \jour Phys. Rev. Lett.
     \vol 56
     \yr 1986
     \pages 889-892
\endref\vskip.08truein

\ref \key Ke1
     \by Kesten, H.
     \pages 125-264
     \paper Aspects of first-passage percolation
     \inbook \'{E}cole d'\'{E}t\'{e} de Probabilit\'{e}s de Saint-Flour
     XIV--1984 
     \ed P. L. Hennequin
     \bookinfo Lecture Notes in Math.
     \vol 1180
     \publ Springer
     \publaddr Berlin-Heidelberg-New York
     \yr 1986 
\endref\vskip.08truein

\ref \key Ke2
     \by Kesten, H.
     \paper On the speed of convergence in first-passage percolation
     \yr 1993
     \jour Ann. Appl. Probab.
     \vol 3
     \pages 296-338
\endref\vskip.08truein

\ref \key KrS
     \by Krug, J., and Spohn, H.
     \paper Kinetic roughening of growing surfaces
     \yr 1991
     \inbook Solids Far from Equilibrium: Growth, Morphology and Defects
     \ed C. Godr\`eche
     \publ Cambridge Univ. Press
     \publaddr Cambridge
\endref\vskip.08truein


\ref \key LN
     \by Licea, C. and Newman, C.M.
     \paper Geodesics in two-dimensional first-passage percolation
     \jour Ann. Probab.  
     \vol 24
     \pages 399-410
     \yr 1996  
\endref\vskip.08truein


\ref \key MR
     \by Meester, R. and Roy, R.
     \book Continuum Percolation
     \publ Cambridge Univ. Press
     \publaddr Cambridge
     \yr 1996
\endref\vskip.08truein

\ref \key N
     \paper Large deviations of sums of independent random variables
     \by Nagaev, S.V.
     \jour Ann. Prob.
     \vol 7
     \yr 1979
     \pages 745-789
\endref\vskip.08truein

\ref \key Ne
     \by Neveu, J. (translated by T. P. Speed)
     \book Martingales a Temps Discret (Discrete-Parameter Martingales)
     \publ Masson \& Cie (American Elsevier)
     \yr 1972 (1975)
     \publaddr Paris (New York)
\endref\vskip.08truein


\ref \key New1
     \by Newman, C.M.
     \paper A Surface View of First-Passage Percolation
     \yr 1995
     \inbook Proceedings of the International Congress of Mathematicians
     \ed S. D. Chatterji
     \publ Birkh\"auser 
     \publaddr Basel-Boston-Berlin
     \pages 1017-1023
\endref\vskip.08truein

\ref \key New2
     \by Newman, C.M.
     \book Topics in Disordered Systems
     \publ Birkh\"auser
     \publaddr Basel-Boston-Berlin
     \yr 1997
\endref\vskip.08truein

\ref \key NewP
     \by Newman, C.M. and Piza, M.S.T.
     \paper Divergence of shape fluctuations in two dimensions
     \jour Ann. Probab. 
     \vol 23
     \pages 977-1005
     \yr 1995
\endref\vskip.08truein

\ref \key NewS1
     \by Newman, C.M. and Stein, D.L.
     \paper Spin glass model with dimension-dependent ground state multiplicity
     \jour Phys. Rev. Lett.
     \vol 72
     \pages 2286-2289
     \yr 1994
\endref\vskip.08truein

\ref \key NewS2
     \by Newman, C.M. and Stein, D.L.
     \paper Ground state structure in a highly disordered spin-glass model
     \jour J. Stat. Phys.
     \vol 92
     \pages 1113-1132
     \yr 1996
\endref\vskip.08truein

\ref \key R
     \by Richardson, D.
     \paper Random growth in a tesselation
     \jour Proc. Cambridge Phil. Soc.
     \yr 1973
     \vol 74
     \pages 515-528 
\endref\vskip.08truein

\ref \key S
     \by Schramm, O.
     \paper Scaling limits of loop-erased random walks and uniform
     spanning trees
     \toappear
     \jour Israel J. Math 
\endref\vskip.08truein

\ref \key Se
     \by Serafini, H.C.
     \book First-passage percolation in the 
     Delaunay graph of a\break $d$-dimensional Poisson process
     \bookinfo Ph.D. Dissertation
     \publ New York University--Courant Inst.~of Math.~Sciences
     \yr 1997
\endref\vskip.08truein

\ref \key SmW
     \by Smythe, R.T. and Wierman, J.C.
     \book First-Passage Percolation on the Square Lattice
     \bookinfo Lecture Notes in Math.
     \vol 671
     \publ Springer
     \publaddr Berlin-Heidelberg-New York
     \yr 1978
\endref\vskip.08truein




\ref \key VW1
     \by Vahidi-Asl, M.Q. and Wierman, J.C.
     \paper First-passage percolation on the Voronoi\break
 tessellation and Delaunay triangulation
     \inbook Random Graphs '87
     \eds M. Karo\'nske, J. Jaworski and A. Ruci\'nski
     \pages 341-359
     \publ Wiley
     \publaddr New York
     \yr 1990
\endref\vskip.08truein


\ref \key VW2
     \by Vahidi-Asl, M.Q. and Wierman, J.C.
     \paper A shape result for first-passage percolation on the Voronoi tessellation and Delaunay triangulation
     \inbook Random Graphs '89
     \eds A. Frieze and T. Luczak
     \pages 247-262
     \publ Wiley
     \publaddr New York
     \yr 1992
\endref\vskip.08truein


\ref \key W
     \by Wehr, J.
     \paper On the number of infinite geodesics and ground states in disordered
     systems
     \jour J. Stat. Phys.
     \vol 87
     \yr 1997
     \pages 439-447
\endref\vskip.08truein

\ref \key Y
     \by Yukich, J.E. 
     \book Probability Theory of Classical Euclidean Optimization Problems
     \bookinfo Lecture Notes in Math. 
     \vol 1675
     \publ Springer
     \publaddr Berlin-Heidelberg-New York
     \yr 1998
\endref\vskip.08truein

\ref \key ZS1
     \by Zuev, S.A. and Sidorenko, A.F.
     \pages 76--86 (51--58 in translation from Russian)
     \paper Continuous models of percolation theory I
     \yr 1985
     \vol 62
     \jour Theoretical and Mathematical Physics
\endref\vskip.08truein    
      
\ref \key ZS2
     \by Zuev, S.A. and Sidorenko, A.F.
     \pages 253--262 (171--177 in translation from Russian)
     \paper Continuous models of percolation theory II
     \yr 1985
     \vol 62
     \jour Theoretical and Mathematical Physics
\endref    




\endRefs
\enddocument
\end